\documentclass{amsart}
\advance\textheight by \topskip

\usepackage{graphicx}
\newtheorem{teo}{\sc Theorem}[section]

\newtheorem{cor}{\sc Corollary}[section]
\newtheorem{lemma}{\sc Lemma}[section]
\newtheorem{pro}{\sc Proposition}[section]

\theoremstyle{definition}
\newtheorem{defi}{\sc Definition}[section]

\theoremstyle{remark}
\newtheorem{rem}[teo]{\sc Remark}

\newtheorem{notation}{\sc Notation}

\newcommand{\bte}{\begin{teo}}
\newcommand{\ete}{\end{teo}}
\newcommand{\bc}{\begin{cor}}
\newcommand{\ec}{\end{cor}}
\newcommand{\bp}{\begin{pro}}
\newcommand{\ep}{\end{pro}}
\newcommand{\bl}{\begin{lemma}}
\newcommand{\el}{\end{lemma}}
\newcommand{\bd}{\begin{defi}}
\newcommand{\ed}{\end{defi}}
\newcommand{\bno}{\begin{notation}}
\newcommand{\eno}{\end{notation}}
\newcommand{\fp}{\hfill $\Box$}

\newcommand{\bca}{\begin{cases}}
\newcommand{\eca}{\end{cases}}
\newcommand{\la}{\langle}
\newcommand{\ra}{\rangle}

\newcommand{\bq}{\begin{equation}}
\newcommand{\eq}{\end{equation}}
\newcommand{\btabu}{\begin{table}}
\newcommand{\etabu}{\end{table}}
\newcommand{\bt}{\begin{tabular}}
\newcommand{\et}{\end{tabular}}
\newcommand{\ba}{\begin{array}}
\newcommand{\ea}{\end{array}}
\newcommand{\br}{\begin{eqnarray}}
\newcommand{\er}{\end{eqnarray}}
\newcommand{\brn}{\begin{eqnarray*}}
\newcommand{\ern}{\end{eqnarray*}}
\newcommand{\benu}{\begin{enumerate}}
\newcommand{\eenu}{\end{enumerate}}
\newcommand{\bite}{\begin{itemize}}
\newcommand{\eite}{\end{itemize}}

\newcommand{\bbc}{\Bbb{C}}

\newcommand{\supp}{\operatorname{supp }}

\title[An introduction to multiple orthogonal polynomials and Hermite-Pad\'e approximation] {An introduction to multiple orthogonal polynomials and Hermite-Pad\'e approximation}

\thanks{This work  was supported by research grant MTM2015-65888-C4-2-P of Ministerio de Econom\'{\i}a  y Competitividad, Spain.}

\author[L\'{o}pez]{G. L\'opez Lagomasino}
\address[L\'{o}pez]{Departamento de Matem\'aticas,
Universidad Carlos III de Madrid, c/ Universidad 30, 28911 Legan\'es, Spain.}

\email[L\'{o}pez]{lago@@math.uc3m.es}

\begin{document}

\maketitle

\begin{abstract}
We present a brief introduction to the theory of multiple orthogonal polynomials on the basis of known results for an important class of measures known as Nikishin systems. For type I and type II multiple orthogonal polynomials with respect to such systems of measures, we describe some of their  most relevant properties regarding location and distribution of zeros as well as their weak and ratio asymptotic behavior.
\end{abstract}

\vspace{1cm}


{\it A.M.S. Subject Classification.} Primary: 30E10, 42C05;
Secondary: 41A20.

\section{Introduction} \label{sec:HP}

The object of this paper is to provide an introduction to the study of Hermite-Pad\'e approximation, multiple orthogonal polynomials, and some of their asymptotic properties. For the most part, the attention is restricted to the case of multiple orthogonality with respect to an important class of measures introduced by E.M. Nikishin in \cite{Nik}. For this reason, they are referred to in the specialized literature as  Nikishin systems.

Throughout the years, with the assistance of colleagues and former Ph. D. students, I have dedicated a great part of my research to this subject area and I wish to express here my gratitude to all those involved. This material is not intended to present new results, but a mere account of the experience I have accumulated for the benefit of future research from the younger generations.

\subsection{Some historical background.}
In 1873, Charles Hermite publishes in \cite{Her} his proof of the transcendence of $e$ making use of simultaneous rational approximation of systems of exponentials. That paper marked the beginning of the modern analytic theory of numbers.

The formal theory of simultaneous rational approximation for general systems of analytic functions was initiated by K. Mahler in lectures delivered at the University of Groningen in 1934-35. These lectures were published years later in \cite{Mah}. Important contributions in this respect are also due to his students J. Coates and H. Jager, see \cite{Coa} and \cite{Jag}. K. Mahler's approach to the simultaneous approximation of finite systems of analytic functions may be reformulated in the following terms.

Let ${\bf f} = (f_1,\ldots,f_m)$ be a  family of analytic functions in some domain $D$ of the extended complex plane containing $\infty$. Fix a non-zero multi-index ${\bf n} = (n_1,\ldots,n_m) \in {\mathbb{Z}}_+^{m}, |{\bf n}| = n_1+\ldots,n_m$. There exist polynomials $a_{{\bf n},1}, \ldots, a_{{\bf n},m}$, not all identically equal to zero, such that
\begin{itemize}
\item[i)] $\deg a_{{\bf n},j} \leq n_j -1, j=1,\ldots,m$ $(\deg a_{{\bf n},j} \leq -1$ means that $a_{{\bf n},j} \equiv 0$),
\item[ii)]  $a_{{\bf n},0} +\sum_{j=1}^m a_{{\bf n},j}(z) f_j(z)  = {\mathcal{O}}(1/z^{|{\bf n}|}), z \to \infty.$
\end{itemize}
for some polynomial $a_{{\bf n},0}$. Analogously, there exists $Q_{\bf n}$, not identically equal to zero, such that
\begin{itemize}
\item[i)] $\deg Q_{{\bf n}} \leq |{\bf n}|,$
\item[ii)] $Q_{\bf n}(z)f_j(z) - P_{{\bf n},j}(z)= {\mathcal{O}}(1/z^{ n_j+1}), z \to \infty, j=1,\ldots,m,$
\end{itemize}
for some polynomials $P_{{\bf n},j}, j=1,\ldots,m.$

The polynomials $a_{{\bf n},0}$ and $P_{{\bf n},j},j=0,\ldots,m,$ are uniquely determined from ii) once their partners are found. The two constructions are called type I and type II polynomials (approximants) of the system $(f_1,\ldots,f_m)$. Algebraically, they are closely related. This is clearly exposed in \cite{Coa}, \cite{Jag}, and \cite{Mah}. When $m=1$ both definitions coincide with that of the well-known Pad\'{e} approximation in its linear presentation.

Apart from Hermite's result, type I, type II, and a combination of the two (called mixed type), have been employed in the proof of the irrationality of other numbers. For example, in \cite{Beu} F. Beukers shows that Apery's  proof (see \cite{Ape}) of the irrationality of $\zeta(3)$  can be placed in the context of mixed type Hermite-Pad\'{e} approximation. See \cite{Ass} for a brief introduction and survey on the subject. More recently, mixed type approximation has appeared in random matrix and non-intersecting brownian motion theories (see, for example, \cite{BleKui}, \cite{DK}, \cite{K1}, and \cite{K2}), and the Degasperi-Procesi equation \cite{Bertola:CBOPs}.

In applications in the areas of number theory, convergence of simultaneous rational approximation, and asymptotic properties of type I and type II polynomials, a central question is if these polynomials have no defect; that is, if they attain the maximal degree possible.

\begin{defi}\label{def1}
A multi-index ${\bf n}$ is said to be {\bf normal} for the system ${\bf f}$ for type I approximation (respectively, for type II,) if $\deg a_{{\bf n},j} = n_j -1, j=1,\ldots,m$ (respectively, $\deg Q_{\bf n} = |{\bf n}|$). A system of functions ${\bf f}$ is said to be {\bf perfect} if all multi-indices are normal.
\end{defi}

It is easy to verify that $(a_{{\bf n},0},\ldots,a_{{\bf n},m})$ and $Q_{\bf n}$ are uniquely determined to within a constant factor when ${\bf n}$ is normal. Moreover, if a system is perfect, the order of approximation in parts ii) above is exact for all ${\bf n}$. The convenience of these properties is quite clear.

Considering the construction at the origin (instead of $z=\infty$ which we chose for convenience), the system of exponentials considered by Hermite, $(e^{w_1z},\ldots,e^{w_mz}), w_i \neq w_j, i \neq j, i,j = 1,\ldots,m, $ is known to be perfect for type I and type II. A second example of a perfect system for both types is that given by the binomial functions $(1-z)^{w_1},\ldots,(1-z)^{w_m}, w_i -w_j \not\in {\mathbb{Z}}$.  When normality occurs for multi-indices with decreasing components the system is said to be weakly perfect. Basically, these are the only examples known of perfect  systems, except for certain ones formed by Cauchy transforms of measures.

\subsection{Markov systems and orthogonality.} Let $s$ be a finite Borel measure with constant sign whose compact support consists of infinitely many points and is contained in the real line.  In the sequel, we only consider such measures. By $\Delta$ we denote the smallest interval which contains the support, $\supp{s},$ of $s$. We denote this class of measures by ${\mathcal{M}}(\Delta)$. Let
\[ \widehat{s}(z) = \int \frac{ds(x)}{z-x}
\]
denote the Cauchy transform of $s$. Obviously, $\widehat{s} \in {\mathcal{H}}(\overline{\mathbb{C}}\setminus \Delta);$  that is, it is analytic in $\overline{\mathbb{C}}\setminus \Delta.$

If we apply the construction above to the system formed by $\widehat{s}\,\, (m=1)$, it is easy to verify that $Q_{\bf n}$ turns out to be orthogonal to all polynomials of degree less than $n \in {\mathbb{Z}}_+$. Consequently, $\deg Q_{\bf n} = n,$ all its zeros are simple and lie in the open convex hull $\mbox{Co}(\supp s)$ of $\supp s$. Therefore, such systems of one function are perfect. These properties allow to deduce Markov's theorem on the convergence of (diagonal) Pad\'e approximations of $\widehat{s}$   published in \cite{Mar}. For this reason, $\widehat{s}$ is also called a Markov function.

Markov functions are quite relevant in several respects. Many elementary functions can be expressed as such. The resolvent function of self-adjoint operators admits that type of representation. If one allows complex weights, any reasonable analytic function in the extended complex plane with a finite number of algebraic singularities adopts that form. This fact, and the use of Pad\'e approximation, has played a central role in some of the most relevant achievements in the last decades concerning the exact rate of convergence of the best rational approximation: namely, A.A. Gonchar and E.A. Rakhmanov's result, see \cite{GR1}, \cite{Gon1}, and \cite{Apt3}, on the best rational approximation of $e^{-x}$ on $[0,+\infty)$; and H. Stahl's theorem, see \cite{S}, on the best rational approximation of $x^{\alpha}$ on $[0,1]$.

Let us see two examples of general systems of Markov functions which play a central role in the theory of multiple orthogonal polynomials.

\subsubsection{Angelesco systems.} In \cite{Ang}, A. Angelesco considered the following systems of functions. Let $\Delta_j, j=1,\ldots,m,$ be pairwise disjoint bounded intervals contained in the real line and $s_j, j=1,\ldots,m,$ a system of measures such that $\mbox{Co}(\supp s_j) = \Delta_j.$

Fix ${\bf n} \in {\mathbb{Z}}_+^{m}$ and consider the type II approximant of the so called Angelesco system of functions $(\widehat{s}_1,\ldots,\widehat{s}_m)$ relative to $\bf n$. It turns out that
\[ \int x^{\nu} Q_{\bf n}(x) ds_j(x) =0, \quad \nu = 0,\ldots,n_j -1,\quad j=1,\ldots,m.
\]
Therefore, $Q_{\bf n}$ has $n_j$ simple zeros in the interior (with respect to the euclidean topology of ${\mathbb{R}}$) of $\Delta_j$. In consequence, since the intervals $\Delta_j$ are pairwise disjoint, $\deg Q_{\bf n} = |{\bf n}|$ and Angelesco systems are type II perfect.

Unfortunately, Angelesco's paper received little attention and such systems  reappear many years later in \cite{Nik1} where E.M. Nikishin  deduces some of their formal properties.

In \cite{GR} and \cite{Apt1}, their logarithmic and strong asymptotic behavior, respectively, are given. These multiple orthogonal polynomials and the rational approximations associated have nice asymptotic formulas but not so good convergence properties. In this respect, a different  system of Markov functions turns out to be more interesting and foundational from the geometric and analytic points of view.

\subsubsection{Nikishin systems.}  In an attempt to construct general classes of functions for which normality takes place, in \cite{Nik} E.M. Nikishin introduced the concept of MT-system (now called Nikishin system). Let $\Delta_{\alpha}, \Delta_{\beta}$ be two non intersecting bounded intervals contained in the real line and $\sigma_{\alpha} \in {\mathcal{M}}(\Delta_{\alpha}), \sigma_{\beta} \in {\mathcal{M}}(\Delta_{\beta})$. With these two measures we define a third one as follows (using the differential notation)
\[ d \la \sigma_{\alpha},\sigma_{\beta} \ra (x) = \widehat{\sigma}_{\beta}(x) d\sigma_{\alpha}(x);
\]
that is, one multiplies the first measure by a weight formed by the Cauchy transform of the second measure. Certainly, this product of measures is non commutative. Above, $\widehat{\sigma}_{\beta}$ denotes the Cauchy transform of the measure $\sigma_{\beta}$.

\begin{defi} \label{Nikishin} Take a collection  $\Delta_j, j=1,\ldots,m,$ of intervals such that
\[ \Delta_j \cap \Delta_{j+1} = \emptyset, \qquad j=1,\ldots,m-1.
\]
Let $(\sigma_1,\ldots,\sigma_m)$ be a system of measures such that $\mbox{Co}(\supp \sigma_j) = \Delta_j, \sigma_j \in {\mathcal{M}}(\Delta_j), j=1,\ldots,m.$
We say that $(s_{1,1},\ldots,s_{1,m}) = {\mathcal{N}}(\sigma_1,\ldots,\sigma_m)$, where
\[ s_{1,1} = \sigma_{1}, \quad s_{1,2} = \la \sigma_1,\sigma_2 \ra, \ldots \quad , s_{1,m} = \la \sigma_1,\la \sigma_2,\ldots,\sigma_m \ra \ra
\]
is the Nikishin system of measures generated by $(\sigma_1,\ldots,\sigma_m)$.
\end{defi}

Fix ${\bf n} \in {\mathbb{Z}}_+^{m}$ and consider the type II approximant of the Nikishin system of functions $(\widehat{s}_{1,1},\ldots,\widehat{s}_{1,m})$ relative to $\bf n$. It is easy to prove that
\[ \int x^{\nu} Q_{\bf n}(x) ds_{1,j}(x) =0, \quad \nu = 0,\ldots,n_j -1,\quad j=1,\ldots,m.
\]
All the measures $s_{1,j}$ have the same support; therefore, it is not immediate to conclude that $\deg Q_{\bf n} = |{\bf n}|$. Nevertheless, if we denote
\[s_{j,k} = \la \sigma_j,\sigma_{j+1},\ldots,\sigma_k \ra,\qquad j < k, \quad s_{j,j} = \la \sigma_j \ra = \sigma_j,
\]
the previous orthogonality relations may be rewritten as follows
\begin{equation}
\label{eq:a}
\int (p_1 (x) + \sum_{k=2}^m p_k(x)\widehat{s}_{2,k}(x))Q_{\bf n}(x) d \sigma_1(x) =0,
\end{equation}
where $p_1,\ldots,p_m$ are arbitrary polynomials such that $\deg p_k \leq n_k -1, k=1,\ldots,m.$

\begin{defi} \label{def:AT} A system of real continuous functions $u_1,\ldots,u_m$ defined on an interval $\Delta$ is called an AT-system on $\Delta$ for the multi-index ${\bf n} \in {\mathbb{Z}}_+^{m}$ if for any choice of real polynomials (that is, with real coefficients) $p_1,\ldots,p_m, \deg p_k \leq n_k -1,$ the function
\[ \sum_{k=1}^m p_k(x) u_k(x)
\]
has at most $|{\bf n}| -1$ zeros on $\Delta$. If this is true for all ${\bf n} \in {\mathbb{Z}}_+^{m}$ we have an AT system on $\Delta$.
\end{defi}

In other words, $u_1,\ldots,u_m$ forms an AT-system for ${\bf n}$ on $\Delta$ when the system of functions
\[ (u_1,\ldots, x^{n_1 -1}u_1,u_2,\ldots,x^{n_m-1}u_m)
\]
is a Tchebyshev system on $\Delta$ of order $|{\bf n}| -1$. From the properties of Tchebyshev systems (see \cite[Theorem 1.1]{KN}), it follows that given
$x_1,\ldots,x_N, N < |{\bf n}|,$ points in the interior of $\Delta$ one can find polynomials $h_1,\ldots,h_m,$ conveniently, with $\deg h_k \leq n_k -1,$ such that $\sum_{k=1}^m h_k(x) u_k(x)$
changes  sign at $x_1,\ldots,x_N,$ and has no other points where it changes sign on $\Delta.$

In \cite{Nik}, Nikishin stated without proof that the system of functions $(1,\widehat{s}_{2,2},\ldots,\widehat{s}_{2,m})$ forms an AT-system for all multi-indices ${\bf n}\in {\mathbb{Z}}_+^{m}$ such that $n_1 \geq n_1 \geq \cdots \geq n_m$ (he proved it when additionally $n_1 -n_m \leq 1$). Due to (\ref{eq:a}) this implies that Nikishin systems are type II weakly perfect. Ever since the appearance of \cite{Nik}, a subject of major interest for those involved in simultaneous approximation was to determine whether or not Nikishin systems are perfect. This problem was settled positively in \cite{FidLop} (see also \cite{FidLop2} where Nikishin systems with unbounded and or touching supports is considered).

In the last two decades, a general theory of multiple orthogonal polynomials and Hermite-Pad\'e approximation has emerged which to a great extent matches what is known to occur for standard orthogonal polynomials and Pad\'e approximation. From the approximation point of view, Markov and Stieltjes type theorems have been obtained (see, for example, \cite{Bus, DrSt1, DrSt2, LF2, LF3, GR, GRS, LS, LS2,Nik}). From the point of view of the asymptotic properties of multiple orthogonal polynomials there are results concerning their weak, ratio, and strong asymptotic behavior (see, for example, \cite{Apt1, Apt2, AptLopRoc,FLLS,LL,LS}). This is specially so for Nikishin systems of measures on which we will focus in this brief introduction.

\section{On the perfectness of  Nikishin systems.}

Let us begin with the following result which was established in \cite{FidLop}. It constitutes the key to the proof of many interesting properties of Nikishin systems; in particular, the fact that they are perfect. From the definition it is obvious that if $(\sigma_1,\ldots,\sigma_m)$ generates a Nikishin system so does $(\sigma_j,\ldots,\sigma_k)$ where $1\leq j < k \leq m$.

\begin{teo}
\label{teo:1}
Let $(s_{1,1},\ldots,s_{1,m}) = {\mathcal{N}}(\sigma_{1},\ldots,\sigma_m)$ be given. Then, the system  $(1,\widehat{s}_{1,1},\ldots,\widehat{s}_{1,m})$ forms an AT-system on any interval $\Delta$ disjoint from $\Delta_1 = \mbox{\rm Co}(\supp \sigma_1)$. Moreover, for each ${\bf n} \in {\mathbb{Z}}_+^{m+1},$ and arbitrary polynomials with real coefficients $p_k, \deg p_k \leq n_k-1, k=0,\ldots,m,$ the linear form $p_0 + \sum_{k=1}^m p_k \widehat{s}_{1,k}, $ has at most $|{\bf n}|-1$ zeros in ${\mathbb{C}} \setminus \Delta_1$.
\end{teo}

For arbitrary multi-indices the proof is quite complicated and based on intricate transformations which allow to reduce the problem to the case of multi-indices with decreasing components. For such multi-indices the proof is pretty straightforward and we will limit ourselves to that situation.  Let us first present two auxiliary lemmas.

\bl \label{reduc}  Let $(s_{1,1},\ldots,s_{1,m}) = \mathcal{N}(\sigma_1,\ldots,\sigma_m)$ be given. Assume that there exist polynomials with real coefficients $\ell_0,\ldots,\ell_m$ and a polynomial $w$ with real coefficients whose zeros lie in $\mathbb{C} \setminus \Delta_1$  such that
\[\frac{\mathcal{L}_0(z)}{w(z)} \in \mathcal{H}(\mathbb{C} \setminus \Delta_1)\qquad \mbox{and} \qquad \frac{\mathcal{L}_0(z)}{w(z)} = \mathcal{O}\left(\frac{1}{z^N}\right), \quad z \to \infty,
\]
where $\mathcal{L}_0  := \ell_0 + \sum_{k=1}^m \ell_k  \widehat{s}_{1,k} $ and $N \geq 1$. Let $\mathcal{L}_1  := \ell_1 + \sum_{k=2}^m \ell_k  \widehat{s}_{2,k} $. Then
\begin{equation} \label{eq:g}
\frac{\mathcal{L}_0(z)}{w(z)} = \int \frac{\mathcal{L}_1(x)}{(z-x)} \frac{{\rm d}\sigma_1(x)}{w(x)}.
\end{equation}
If $N \geq 2$, we also have
\begin{equation} \label{eq:f}
\int x^{\nu}  \mathcal{L}_1(x)  \frac{{\rm d}\sigma_1(x)}{w(x)} = 0, \qquad \nu = 0,\ldots, N -2.
\end{equation}
In particular, $\mathcal{L}_1$ has at least $N -1$ sign changes in  the interior of  $ {\Delta}_1 $ (with respect to the Euclidean topology of $\mathbb{R}$).
\el

{\bf Proof.} Let $\Gamma $ be a positively oriented closed
smooth Jordan curve that surrounds $\Delta_1$ sufficiently close to $\Delta_1$.
Since $\frac{\mathcal{L}_0(z)}{w(z)}  =
{\mathcal{O}}(1/z), z \to \infty,$  if $z$ and the zeros of $w(z)$ are in the
unbounded connected component of the complement of $\Gamma$, Cauchy's
integral formula and Fubini's theorem render
\[ \frac{\mathcal{L}_0(z)}{w(z)}  = \frac{1}{2\pi i} \int_{\Gamma} \frac{\mathcal{L}_0(\zeta)}{w(\zeta)} \frac{d\zeta}{z-\zeta} = \frac{1}{2\pi i} \sum_{k=1}^m \int_{\Gamma}
\frac{\ell_k(\zeta)\widehat{s}_{1,k}(\zeta)  d\zeta}{w(\zeta)(z-\zeta)} = \] \[
\sum_{k=1}^m \int \frac{1}{2\pi i} \int_{\Gamma} \frac{\ell_k(\zeta)
d\zeta}{w(\zeta)(z-\zeta)(\zeta -x)} ds_{1,k}(x) = \sum_{k=1}^m \int
\frac{\ell_k(x) ds_{1,k}(x)}{w(x)(z-x)} = \int \frac{\mathcal{L}_1(x)}{(z-x)} \frac{{\rm d}\sigma_1(x)}{w(x)}
\]
which is (\ref{eq:g}).

When $N \geq 2$, it follows that $\frac{z^\nu\mathcal{L}_0(z)}{w(z)}  =
{\mathcal{O}}(1/z^2), z \to \infty,$ for $\nu =0,\ldots,N-2$. Then, using Cauchy's theorem, Fubini's
theorem and Cauchy's integral formula, it follows that
\[ 0 = \int_{\Gamma}\frac{z^\nu\mathcal{L}_0(z)}{w(z)}  dz = \sum_{k=1}^m \int_{\Gamma} \frac{z^\nu\ell_k(z)\widehat{ s}_{1,k}(z)}{w(z)} dz = \sum_{k=1}^m \int
\int_{\Gamma} \frac{z^\nu \ell_k(z)dz}{w_k(z)(z-x)}ds_{1,k}(x) =
\]
\[ 2\pi i \sum_{k=1}^m \int
\frac{x^\nu\ell_k(x) }{w_k(x) }  ds_{1,k}(x) = 2\pi i \int x^{\nu}  \mathcal{L}_1(x)  \frac{{\rm d}\sigma_1(x)}{w(x)} ,
\]
and we obtain (\ref{eq:f}).
\fp

\bl
\label{lem:3}
Let $(s_{1,1},\ldots,s_{1,m}) =
{\mathcal{N}}(\sigma_1,\ldots,\sigma_m)$ and ${\bf n} = (n_0,\ldots,n_m) \in
{\mathbb{Z}}_+^{m+1}$ be given. Consider the linear form
\[ {\mathcal{L}}_{\bf n}  = p_0 + \sum_{k=1}^m p_k \widehat{s}_{1,k}, \quad \deg p_k \leq n_k-1, \quad k=0,\ldots,m,
\]
where the polynomials $p_k$ have real coefficients. Assume that $n_0 = \max\{n_0,n_1-1,\ldots,n_m-1\}$. If
${\mathcal{L}}_{\bf n}$ had at least $|{\bf n}|$ zeros in $ {\mathbb{C}} \setminus \Delta_1$ the reduced form
$p_1 + \sum_{k=2}^m p_k \widehat{s}_{2,k}$ would have at least $|{\bf n}| - n_0$ zeros in $ {\mathbb{C}} \setminus \Delta_2$.
\el

{\bf Proof.} The function ${\mathcal{L}}_{\bf n}$ is symmetric with respect to the real line ${\mathcal{L}}_{\bf n}(\overline{z}) = \overline{{\mathcal{L}}_{\bf n}(z)}$; therefore, its zeros come in conjugate pairs. Thus, if ${\mathcal{L}}_{\bf n}$ has at least $|{\bf n}|$ zeros in $ {\mathbb{C}} \setminus \Delta_1$, there exists a polynomial $w_{\bf n}, \deg w_{\bf n}
\geq |{\bf n}|,$  with real coefficients and zeros contained in
${\mathbb{C}}\setminus \Delta_1$ such that ${\mathcal{L}}_{\bf
n}/w_{\bf n} \in {\mathcal{H}}({\mathbb{C}}\setminus \Delta_1)$. This function has a zero of order $ \geq |{\bf n}| - n_0 +1$
at $\infty$. Consequently, for all $ \nu
=0,\ldots,|{\bf n}| - n_0 -1,$
\[ \frac{z^{\nu} {\mathcal{L}}_{\bf n}}{w_{\bf n}} =
{\mathcal{O}}(1/z^2) \in {\mathcal{H}}(\overline{\mathbb{C}} \setminus
\Delta_1), \qquad z \to \infty, \qquad
\]
and
\[ \frac{z^{\nu} {\mathcal{L}}_{\bf n}}{w_{\bf n}} =  \frac{z^{\nu}p_0}{w_{\bf
n}} +  \sum_{k=1}^m \frac{z^{\nu}  p_k }{w_{\bf n}}
\widehat{s}_{1,k}\,.
\]
From (\ref{eq:f}), it follows that
\[ 0 = \int x^{\nu}  ( p_1 + \sum_{k=2}^m p_k
\widehat{s}_{2,k})(x) \frac{d \sigma_1(x)}{w_{\bf n}(x)}, \qquad \nu
=0,\ldots,|{\bf n}| - n_0 -1,
\]
taking into consideration that $s_{1,1} = \sigma_1$ and $ds_{1,k}(x) =
\widehat{s}_{2,k}(x) d\sigma_1(x), k=2,\ldots,m$.

These orthogonality relations imply that $p_1 + \sum_{k=2}^m p_k
\widehat{s}_{2,k}$ has at least $|{\bf n}| - n_0$ sign changes in the interior of $\Delta_1$.  In fact, if there were at most $|{\bf n}| - n_0 -1$ sign changes one can easily construct a polynomial $p$ of degree $\leq |{\bf n}| - n_0 -1$ such that $p(p_1 + \sum_{k=2}^m p_k
\widehat{s}_{2,k})$ does not change sign on $\Delta_1$ which contradicts the orthogonality relations. Therefore, already in the interior of $\Delta_1 \subset {\mathbb{C}} \setminus \Delta_2$, the reduced form would have the number of zeros claimed. \hfill $\Box$ \vspace{0,2cm}

{\bf Proof of Theorem \ref{teo:1} when $n_0 \geq n_1\geq \cdots \geq n_m$.} In this situation assume that the linear form $p_0 + \sum_{k=1}^m p_k \widehat{s}_{1,k}, $ has at least $|{\bf n}|$ zeros in $\overline{\mathbb{C}} \setminus \Delta_1$. Applying Lemma \ref{lem:3} consecutively $m$ times  we would arrive to the conclusion that $p_m$ has at least $n_m$ zeros in $\mathbb{C}$ but this is impossible since its degree is $\leq n_m -1$. \hfill $\Box$

\medskip

From Theorem \ref{teo:1} the following result readily follows.

\begin{teo}\label{teo:2} Nikishin systems are type I and type II perfect.
\end{teo}

{\bf Proof.} Consider a Nikishin system ${\mathcal{N}}(\sigma_1,\ldots,\sigma_m)$. Let us prove that it is  type I perfect. Given a multi-index $\bf n \in {\mathbb{Z}}_+^{m}$ condition ii) for the system of functions $(\widehat{s}_{1,1},\ldots,\widehat{s}_{1,m})$ implies that
\[ \int x^{\nu} \left(a_{{\bf n},1} + a_{{\bf n},2} \widehat{s}_{2,2} +\cdots + a_{{\bf n},m} \widehat{s}_{2,m}\right)(x) d\sigma_1(x) = 0, \qquad \nu = 0,\ldots, |{\bf n}| -2.
\]
These orthogonality relations imply that ${\mathcal{A}}_{{\bf n},1} := a_{{\bf n},1} + a_{{\bf n},2} \widehat{s}_{2,2} +\cdots + a_{{\bf n},m} \widehat{s}_{2,m}$ has at least $|{\bf n}| -1$ sign changes in the interior of $\Delta_1$ (with the Euclidean topology of $\mathbb{R}$). Suppose that  ${\bf n}$ is not normal; that is, $\deg a_{{\bf n},j}\leq n_j -2$ for some $j =1,\ldots,m$. Then, according to Theorem \ref{teo:1}, ${\mathcal{A}}_{{\bf n},1}$ can have at most $|{\bf n}| -2$ zeros in ${\mathbb{C}} \setminus \Delta_2$. This contradicts the previous assertion, so such a component $j$ cannot exist in $\bf n$.

The proof of type II perfectness is analogous. Suppose that there exists an ${\bf n}$ such that $Q_{\bf n}$ has less than $|{\bf n}|$ sign changes in the interior of $\Delta_1$. Let $x_k, k =1,\ldots x_N, N \leq |{\bf n}| -1$ be the points where it changes sign. Construct a linear form
\[ p_{1} + p_{ 2} \widehat{s}_{2,2} +\cdots + p_{m} \widehat{s}_{2,m}, \qquad \deg p_k \leq n_k -1,\qquad k=1,\ldots,m.
\]
with a simple zero at each of the points $x_k$ and a zero of multiplicity $|{\bf n}| - N -1$ at one of the end points of $\Delta_1$. This is possible because there are sufficient free parameters in the coefficients of the polynomials $p_k$ and the form is analytic on  a neighborhood of $\Delta_1$. By Theorem \ref{teo:1} this linear form cannot have any more zeros in the complement of $\Delta_2$ than those that have been assigned. However, using
\eqref{eq:a}, we have
\[\int (p_1 (x) + \sum_{k=2}^m p_k(x)\widehat{s}_{2,k}(x))Q_{\bf n}(x) d \sigma_1(x) =0,
\]
which is not possible since the function under the integral sign has constant sign on $\Delta_1$ and it is not identically equal to zero. \hfill $\Box$

\medskip

\section{On the interlacing property of zeros.}

In the sequel, we restrict our attention to multi-indices in
\[  {\mathbb{Z}}_+^{m}(\bullet) := \{{\bf n} \in {\mathbb{Z}}_+^{m}:  n_1 \geq \cdots \geq n_m\}.
\]

\subsection{Interlacing for type I}
Fix ${\bf n} \in {\mathbb{Z}}_+^{m}(\bullet)$. Consider the type I Hermite-Pad\'e approximant $(a_{{\bf n},0},\ldots,a_{{\bf n},m})$ of $(\widehat{s}_{1,1},\ldots,\widehat{s}_{1,m})$ for the  multi-index  ${\bf n}= (n_1,\ldots,n_m)$. Set
\[ {\mathcal{A}}_{{\bf n},k} = a_{{\bf n},k} + \sum_{j=k+1}^m a_{{\bf n},j} \widehat{s}_{k+1,j}, \qquad k=0,\ldots,m-1.
\]
We take ${\mathcal{A}}_{{\bf n},m} = a_{{\bf n},m} $.

\begin{pro} \label{teo:3} For each $k=1,\ldots,m$ the linear form ${\mathcal{A}}_{{\bf n},k}$ has exactly $n_k + \cdots +n_m -1$ zeros in ${\mathbb{C}} \setminus \Delta_{k+1}$, where $\Delta_{m+1}= \emptyset$, they are all simple and lie in the interior of $\Delta_k$. Let $A_{{\bf n},k}$ be the monic polynomial whose roots are the zeros  of ${\mathcal{A}}_{{\bf n},k}$ on $\Delta_k$. Then
\begin{equation} \label{orto1}
\int x^\nu {\mathcal{A}}_{{\bf n},k}(x) \frac{d \sigma_k (x)}{A_{{\bf n},k-1}(x)} = 0, \qquad \nu =  n_k+\cdots+n_m -2, \qquad k=1,\ldots,m,
\end{equation}
where $A_{{\bf n},0} \equiv 1$.
\end{pro}

{\bf Proof.} According to Theorem \ref{teo:1} applied to the Nikishin system ${\mathcal{N}}(\sigma_{k+1},\ldots,\sigma_m)$, the linear form ${\mathcal{A}}_{{\bf n},k}$ cannot have more than $n_k +\cdots,n_m -1$ zeros in ${\mathbb{C}} \setminus \Delta_{k+1}$. So it suffices to show that it has at least $n_k + \cdots + n_m -1$ sign changes on $\Delta_k$ to prove that it has exactly $n_k + \cdots + n_m -1$ simple roots in ${\mathbb{C}} \setminus \Delta_{k+1}$ which lie in the interior of $\Delta_k$. We do this producing consecutively the orthogonality relations \eqref{orto1}.

From the definition of type I Hermite-Pad\'e approximation it follows that $z^{\nu}{\mathcal{A}}_{{\bf n},0} = {\mathcal{O}}(1/z^2), \nu = 0,\ldots,|{\bf n}| -2$. From \eqref{eq:f} we get that
\[\int x^\nu {\mathcal{A}}_{{\bf n},1}(x) d \sigma_1 (x)  = 0, \qquad \nu =  n_1+\cdots+n_m -2,
\]
which is \eqref{orto1} when $k=1$. Therefore, ${\mathcal{A}}_{{\bf n},1}$ has at least $n_1+\cdots+n_m -1$ sign changes on $\Delta_1$ as we needed to prove. Let $A_{{\bf n},1}$ be the monic polynomial whose roots are the zeros of ${\mathcal{A}}_{{\bf n},1}$ on $\Delta_1$. Since $z^\nu {\mathcal{A}}_{{\bf n},1}/A_{{\bf n},1} = {\mathcal{O}}(1/z^2), \nu = n_2+\cdots+n_m -2$ (recall that the multi-index has decreasing components) from \eqref{eq:f} we get \eqref{orto1} for $k=2$ which implies that ${\mathcal{A}}_{{\bf n},2}$ has at least $n_2+\cdots+n_m -1$ sign changes on $\Delta_2$ as needed. We repeat the process until we arrive to ${\mathcal{A}}_{{\bf n},m} =a_{{\bf n},m} = A_{{\bf n},m}$. \hfill $\Box$

\medskip

Fix $\ell \in \{1,\ldots,m\}$. Given ${\bf n} \in {\mathbb{Z}}_+^{m}(\bullet)$, we denote
\[{\bf n}^{\ell} = {\bf n} + {\bf e}_\ell,
\]
where ${\bf e}_\ell$ is the $m$ dimensional canonical vector with $1$ in its $\ell$-th component and $0$ everywhere else.
the multi-index which is obtained adding $1$ to the $\ell$-th component of ${\bf n}$.  Notice that ${\bf n}^{\ell}$ need not belong to ${\mathbb{Z}}_+^{m}(\bullet)$; however, we can construct the linear forms ${\mathcal{A}}_{{\bf n}^{\ell},k}$ corresponding to the type I Hermite-Pad\'e approximation with respect to  ${\bf n}^\ell$.

\begin{teo} \label{teo:4} Assume that ${\bf n} \in {\mathbb{Z}}_+^{m}(\bullet)$ and $\ell \in \{1,\ldots,m\}$. For each $k=1,\ldots,m$,   ${\mathcal{A}}_{{\bf n}^\ell,k}$ has at most   $n_k+\cdots+n_m$ zeros in  ${\mathbb{C}} \setminus \Delta_{k+1}$ and at least $n_k+\cdots+n_m-1$ sign changes in $\Delta_k$. Therefore, all its zeros are real and simple. The zeros of ${\mathcal{A}}_{{\bf n},k}$ and ${\mathcal{A}}_{{\bf n}^\ell,k}$ in  ${\mathbb{C}} \setminus \Delta_{k+1}$ interlace.
\end{teo}
{\bf Proof.} Let $A,B$ be real constants such that $|A|+|B| > 0$. For $k=1,\ldots,m$, consider the forms $G_{{\bf n},k}:= A {\mathcal{A}}_{{\bf n},k} + B{\mathcal{A}}_{{\bf n}^\ell,k}$. From Theorem \ref{teo:1} it follows that $G_{{\bf n},k}$ has at most $n_k+\cdots+n_m$ zeros in  ${\mathbb{C}} \setminus \Delta_{k+1}$. Let us prove that it has at least $n_k+\cdots+n_m -1$ sign changes on $\Delta_k$. Once this is achieved we know that all the zeros of $G_{{\bf n},k}$ are simple and lie on the real line. In particular, this would be true for ${\mathcal{A}}_{{\bf n}^\ell,k}$.

Let us start with $k=1$. Notice that
\[\int x^\nu G_{{\bf n},1}(x) d \sigma_1 (x)  = 0, \qquad \nu =  n_1+\cdots+n_m -2.
\]
Consequently, $G_{{\bf n},1}$ has at least $n_1+\cdots+n_m -1$ sign changes on $\Delta_1$ as claimed. Therefore, its zeros in $\mathbb{C} \setminus \Delta_2$ are real and simple. Let $W_{{\bf n},1}$ be the monic polynomial whose roots are the simple zeros of $G_{{\bf n},1}$ in $\mathbb{C} \setminus \Delta_2$.

Observe that $z^{\nu} G_{{\bf n},1}/W_{{\bf n},1} ={\mathcal{O}}(1/z^2), \nu = 0,\ldots,n_2+\cdots+n_m -2$. Consequently, from \eqref{eq:f}
\[
\int x^\nu G_{{\bf n},2}(x) \frac{d \sigma_2 (x)}{W_{{\bf n},1}(x)} = 0, \qquad \nu =  n_2+\cdots+n_m -2,
\]
which implies that $G_{{\bf n},2}$ has at least $n_2+\cdots+n_m -1$ sign changes on $\Delta_2$. So all its zeros are real and simple. Repeating the same arguments we obtain the claim about the number of sign changes of $G_{{\bf n},k}$ on $\Delta_k$ and that all its zeros are real and simple for each $k=1,\ldots,m$.

Let us check that ${\mathcal{A}}_{{\bf n},k}$ and ${\mathcal{A}}_{{\bf n}^\ell,k}$ do not have common zeros in ${\mathbb{C}}\setminus \Delta_{k+1}$. To the contrary, assume that $x_0$ is such a common zero. We have that ${\mathcal{A}}_{{\bf n},k}'(x_0) \neq 0 \neq {\mathcal{A}}_{{\bf n}^\ell,k}'(x_0)$ because the zeros of these forms are simple. Then
\[{\mathcal{A}}_{{\bf n}^\ell,k}'(x_0){\mathcal{A}}_{{\bf n},k} -  {\mathcal{A}}_{{\bf n},k}'(x_0){\mathcal{A}}_{{\bf n}^\ell,k}\]
has a double zero at $x_0$ against what was proved above.

Fix $y \in {\mathbb{R}} \setminus \Delta_{k+1}$ and set
\[ G^y_{{\bf n},k}(z) = {\mathcal{A}}_{{\bf n}^l,k}(z){\mathcal{A}}_{{\bf n},k}(y)-{\mathcal{A}}_{{\bf
n}^l,k}(y){\mathcal{A}}_{{\bf n},k}(z).
\]
Let $x_{1}, x_{2}, x_1 < x_2,$ be two consecutive zeros of ${\mathcal{A}}_{{\bf n}^l,k}$ in
${\mathbb{R}} \setminus \Delta_{k+1}$ and let $y \in
(x_1,x_2)$.
The function $G^y_{{\bf n},k}(z)$ is a real valued function when
restricted to ${\mathbb{R}} \setminus \Delta_{k+1}$ and analytic in
${\mathbb{C}} \setminus \Delta_{k+1}$. We have  $(G^y_{{\bf
n},k})^{\prime}(z) = {\mathcal{A}}_{{\bf n}^\ell,k}^{\prime}(z){\mathcal{A}}_{{\bf
n},k}(y)- {\mathcal{A}}_{{\bf n}^\ell,k}(y){\mathcal{A}}_{{\bf n},k}^{\prime}(z)$. Assume that $(G^{y_0}_{{\bf n},k})^{\prime}(y_0) = 0$ for some $y_0 \in (x_1,x_2)$. Since
$G^{y}_{{\bf n},k}(y) = 0$ for all $y \in (x_1,x_2)$
we obtain that $G^{y_0}_{{\bf n},k}(z)$ has a zero of order
$\geq 2$ (with respect to $z$) at $y_0$ which contradicts what was proved above.  Consequently,
\[(G^{y}_{{\bf n},k})^{\prime}(y) =
{\mathcal{A}}_{{\bf n}^l,k}^{\prime}(y){\mathcal{A}}_{{\bf n},k}(y)- {\mathcal{A}}_{{\bf
n}^l,k}(y){\mathcal{A}}_{{\bf n},k}^{\prime}(y)
\]
takes values with constant sign for all $y \in
(x_1,x_2)$. At the end points $x_1, x_2$,
this function cannot be equal to zero because
${\mathcal{A}}_{{\bf n},k}, {\mathcal{A}}_{{\bf n}^l,k}$ do not have common zeros.
By continuity, $(G^{y}_{{\bf n},k})^{\prime}$ preserves the same
sign on all $[x_{\nu},x_{\nu+1}]$ (and, consequently, on each side
of the interval $\Delta_{k+1}$). Thus
\[ \mbox{sign}(G^{x_1}_{{\bf n},k})^{\prime}(x_1) =
\mbox{sign}(({\mathcal{A}}_{{\bf n}^l,k})^{\prime}(x_1){\mathcal{A}}_{{\bf
n},k}(x_1)) =
\]
\[ \mbox{sign}(({\mathcal{A}}_{{\bf
n}^l,k})^{\prime}(x_2){\mathcal{A}}_{{\bf n},k}(x_2)) =
\mbox{sign}(G^{x_2}_{{\bf n},k})^{\prime}(x_2)\,.
\]
Since
\[ \mbox{sign}({\mathcal{A}}_{{\bf n}^l,k})^{\prime}(x_1) \neq \mbox{sign}({\mathcal{A}}_{{\bf
n}^l,k})^{\prime}(x_2)\,,
\]
we obtain that
\[ \mbox{sign}{\mathcal{A}}_{{\bf n},k}(x_1) \neq \mbox{sign}{\mathcal{A}}_{{\bf
n},k} (x_2)\,;
\]
consequently, there must be an intermediate zero of ${\mathcal{A}}_{{\bf
n},k}$ between $x_1$ and $x_2$.
\hfill $\square$

\medskip

The linear forms  ${\mathcal{A}}_{{\bf n},k}, k=0,\ldots,m$ satisfy nice iterative integral representations. Let us introduce the following functions
\begin{equation} \label{Hk}
H_{{\bf n},k}(z) := \frac{A_{{\bf n},k+1}(z){\mathcal{A}}_{{\bf n},k}(z)}{A_{{\bf n},k}(z)}, \qquad k=0,\ldots,m,
\end{equation}
where the $A_{{\bf n},k}$ are the polynomials introduced in the statement of Proposition \ref{teo:3}. We take $A_{{\bf n},0}\equiv 1 \equiv A_{{\bf n},m+1}$. We have:

\begin{pro} \label{teo:5} For each $k=0,\ldots,m-1$
\begin{equation} \label{orto2}
H_{{\bf n},k}(z) = \int \frac{A_{{\bf n},k+1}^2(x)}{z-x}  \frac{H_{{\bf n},k+1}(x)d \sigma_{k+1} (x)}{A_{{\bf n},k}(x)A_{{\bf n},k+2}(x)}
\end{equation}
and
\begin{equation} \label{orto1*}
\int x^\nu A_{{\bf n},k+1}(x)  \frac{H_{{\bf n},k+1}(x)d \sigma_{k+1} (x)}{A_{{\bf n},k}(x)A_{{\bf n},k+2}(x)} = 0, \qquad \nu =  n_{k+1}+\cdots+n_m -2.
\end{equation}
\end{pro}

{\bf Proof.} Relation \eqref{orto1*} is \eqref{orto1} (with the index $k$ shifted by one) using the notation introduced for the functions $H_{{\bf n},k}(z)$. From \eqref{eq:g} applied to the function ${\mathcal{A}}_{{\bf n},k}/A_{{\bf n},k}, k=0,\ldots,m-1$, we obtain
\[ \frac{{\mathcal{A}}_{{\bf n},k}(z)}{A_{{\bf n},k}(z)} = \int \frac{{\mathcal{A}}_{{\bf n},k+1}(x)}{z-x} \frac{d\sigma_{k+1}(x)}{A_{{\bf n},k}(x)}.
\]
Using \eqref{orto1} with $k+1$ replacing $k$ it follows that
\[ \int \frac{A_{{\bf n},k+1}(z) - A_{{\bf n},k+1}(x)}{z-x} {\mathcal{A}}_{{\bf n},k+1}(x) \frac{d \sigma_{k+1} (x)}{A_{{\bf n},k}(x)} = 0.
\]
Combining these two integral formulas we get
\[ \frac{A_{{\bf n},k+1}(z){\mathcal{A}}_{{\bf n},k}(z)}{A_{{\bf n},k}(z)} = \int \frac{A_{{\bf n},k+1}(x)}{z-x} {\mathcal{A}}_{{\bf n},k+1}(x) \frac{d \sigma_{k+1} (x)}{A_{{\bf n},k}(x)} =
\]
\[\int \frac{A_{{\bf n},k+1}^2(x)}{z-x} \frac{A_{{\bf n},k+2}(x){\mathcal{A}}_{{\bf n},k+1}(x)}{A_{{\bf n},k+1}(x)} \frac{d \sigma_{k+1} (x)}{A_{{\bf n},k}(x)A_{{\bf n},k+2}(x)}
\]
which is \eqref{orto2}. \hfill $\Box$

\medskip

Notice that the varying measure appearing in \eqref{orto2} has constant sign on $\Delta_{k+1}$.

\medskip

\subsection{Interlacing for type II} Now, let us see what happens with the type II Hermite Pad\'e approximants. Let us introduce the following functions.
\[
\Psi_{{\bf n},0}(z)=Q_{{\bf n}}(z),\qquad\Psi_{{\bf n},k}(z)=\int
\frac{\Psi_{{\bf n},k-1}(x)}{z-x} d\sigma_k(x),\qquad k=1,\ldots,m.
\]
where $Q_{\bf n}$ is the type II multiple orthogonal polynomial. The next result is  \cite[Proposition 1]{GRS}.

\begin{pro} \label{teo:6} Let ${\bf n} \in \mathbb{Z}_+^m(\bullet)$. For each $k=1,\ldots,m$
\begin{equation}               \label{eq:fiortogonal}
\int \Psi_{{\bf n},k-1}(x)\left(p_k(x) + \sum_{j=k+1}^m p_j(x) \widehat{s}_{k+1,j}(x)\right) d \sigma_k(x)  =0,\qquad
\deg p_j \leq n_j-1 \qquad j=k,\ldots,m.
\end{equation}
For $k=0,\ldots,m-1$, the function $\Psi_{{\bf n},k}$ has exactly $n_{k+1}+\cdots+n_m$ zeros in
$\mathbb{C} \setminus \Delta_{ k}$, where $\Delta_0 =\emptyset)$, they are all simple and lie in the interior of $\Delta_{k+1}$. The function $\Psi_{{\bf n},m}$ has no roots in $\mathbb{C} \setminus \Delta_{ m}$.
\end{pro}

{\bf Proof.}  The statement about the zeros of $ Q_{\bf n} = \Psi_{{\bf n},0}$ was proved above. The proof of \eqref{eq:fiortogonal} is carried out by induction on $k$. The statement is equivalent to showing that for each $k=1,\ldots,m$, we have
\begin{equation} \label{orto3} \int x^\nu \Psi_{{\bf n},k-1}(x)  d s_{k,j}(x)  =0,\qquad  \nu=0,\ldots,n_j -1, \qquad j=k,\ldots,m.
\end{equation}
When $k=1$ we have $\Psi_{{\bf n},0} = Q_{\bf n}$ and \eqref{orto3} reduces to the orthogonality relations which define $Q_{\bf n}$.
So the basis of induction is settled.

Assume that \eqref{orto3} is true for $k \leq m-1$ and let us show that it also holds for $k+1$. Take $k+1\leq j \leq m$ and $\nu \leq n_j -1$. Then
\[ \int x^\nu \Psi_{{\bf n},k}(x)  d s_{k+1,j}(x)  = \int x^\nu \int \frac{\Psi_{{\bf n},k-1}(t)d\sigma_k(t)}{x-t}  d s_{k+1,j}(x) = \]
\[ \int \Psi_{{\bf n},k-1}(t)\int \frac{ x^\nu - t^\nu + t^\nu}{x-t}  d s_{k+1,j}(x) d\sigma_k(t) =\]
\[  \int \Psi_{{\bf n},k-1}(t) p_{\nu -1}(t)  d\sigma_k(t)  + \int t^\nu \Psi_{{\bf n},k-1}(t) \widehat{s}_{k+1,j}(t) d\sigma_k(t),
\]
where $p_{\nu-1}$ is a polynomial of degree $\leq \nu-1 \leq n_j-2 < n_k -1$ (because the indices have decreasing components).
Since $\widehat{s}_{k+1,j}(t) d\sigma_k(t) = d s_{k,j}(t)$ the induction hypothesis renders that both integrales on the last line equal zero and we obtain what we need.

Applying \eqref{orto3} with $j=k$, we have
\[ \int \frac{z^{n_k} - x^{n_k}}{z-x}  \Psi_{{\bf n},k-1}(x) d\sigma_k(x) = 0.\]
Therefore, using the definition of $\Psi_{{\bf n},k}$, it follows that
\[ z^{n_k}  \Psi_{{\bf n},k}(z) = \int \frac{x^{n_k}}{z-x}  \Psi_{{\bf n},k-1}(x) d\sigma_k(x) = \mathcal{O}(1/z), \qquad z \to \infty.
\]
In other words
\[ \Psi_{{\bf n},k}(z) = \mathcal{O}(1/z^{n_k+1}), \qquad z \to \infty.
\]

Relations \eqref{eq:fiortogonal} applied with $k$ replaced with $k+1\leq m$, together with Theorem \ref{teo:1}, imply that $\Psi_{{\bf n},k}$ has at least $n_{k+1}+\cdots+n_m$ sign changes on the interval $\Delta_{k+1}$. Let $Q_{{\bf n},k+1}$ be the monic polynomial whose roots are the zeros of $\Psi_{{\bf n},k}$ in $\mathbb{C} \setminus \Delta_k$. Take $Q_{{\bf n},m+1} \equiv 1$. Obviously $\deg Q_{{\bf n},k+1} \geq N_{k+1} = n_{k+1} + \cdots, n_m, N_{m+1} = 0$.

Notice that $\Psi_{{\bf n},k}/Q_{{\bf n},k+1} \in \mathcal{H}(\mathbb{C} \setminus \Delta_{k})$ and $\Psi_{{\bf n},k}/Q_{{\bf n},k+1} = \mathcal{O}(1/z^{n_k+N_{k+1}+1}), z \to \infty,$ and the order is $>N_k+ 1$ if $\deg Q_{{\bf n},k+1} > N_{k+1}$. Similar to the way in which \eqref{eq:g} and \eqref{eq:f} were proved, it follows that
\begin{equation} \label{rel1} \frac{\Psi_{{\bf n},k}(z)}{Q_{{\bf n},k+1}(z)} = \int  \frac{\Psi_{{\bf n},k-1}(x)}{z-x}\frac{d\sigma_k(x)}{Q_{{\bf n},k+1}(x)}
\end{equation}
and
\begin{equation} \label{rel2} \int  x^\nu  {\Psi_{{\bf n},k-1}(x)} \frac{d\sigma_k(x)}{Q_{{\bf n},k+1}(x)} = 0, \qquad \nu = 0,\ldots,N_k -1.
\end{equation}
The second one of these relations also implies that $\Psi_{{\bf n},k-1}$ has at least $N_k$ sign changes on $\Delta_k$. Notice that should $\Psi_{{\bf n},k}$ have more than $N_{k+1}$ zeros in $\mathbb{C} \setminus \Delta_{ k}$ we get at least one more orthogonality relation in \eqref{rel2} and we would have that $\Psi_{{\bf n},k-1}$ has more than $N_k$ zeros in $\mathbb{C} \setminus \Delta_{ k-1}, \Delta_0 = \emptyset$. Applying this argument for decreasing values of $k$ we obtain that if for some $k = 1,\ldots,m$, $\Psi_{{\bf n},k}$ has more than $N_{k+1}$ zeros in  $\mathbb{C} \setminus \Delta_{ k}$ then $\Psi_{{\bf n},0} = Q_{\bf n}$ would have more that $|{\bf n}| = N_1 = n_1 +\cdots,n_m$ zeros in $\mathbb{C}$ which is not possible. Consequently, for $k=1,\ldots,m$ the function $\Psi_{{\bf n},k}$ has exactly $N_{k+1}$ zeros in
$\mathbb{C} \setminus \Delta_{ k}$ they are all simple and lie in the interior of $\Delta_{k+1}, \Delta_{m+1} = \emptyset$ as stated. \hfill $\Box$

\medskip

Set
\[{\mathcal{H}}_{{\bf n},k} :=  \frac{Q_{{\bf n},k-1}\Psi_{{\bf n},k-1}}{Q_{{\bf
n},k}}\,.\]

\begin{pro} \label{pro1}
Fix ${\bf n} \in \mathbb{Z}_+^m(\bullet)$. For each $k=1,\ldots,m$
\begin{equation} \label{rel3}
\int x^{\nu} Q_{{\bf n},k}(x)\frac{\mathcal{H}_{{\bf
n},k}(x)d\sigma_k(x)}{Q_{{\bf n},k-1}(x) Q_{{\bf
n},k+1}(x)}=0, \qquad \nu = 0,\ldots,n_k +\cdots+n_m -1,
\end{equation}
and
\begin{equation} \label{rel4}
\mathcal{H}_{{\bf n},k+1}(z) =  \int \frac{Q_{{\bf n},k}^2(x)}{z-x}\frac{\mathcal{H}_{{\bf
n},k}(x)d\sigma_k(x)}{Q_{{\bf n},k-1}(x) Q_{{\bf
n},k+1}(x)},
\end{equation}
where $Q_{{\bf n},0} = Q_{{\bf n},m+1} \equiv 1$.
\end{pro}

{\bf Proof.} Using the notation introduced for the functions $\mathcal{H}_{{\bf n},k}$, \eqref{rel2} adopts the form \eqref{rel3}.
In turn, this implies that
\[ \int \frac{Q_{{\bf n},k}(z)- Q_{{\bf n},k}(x)}{z-x} Q_{{\bf n},k}(x)\frac{\mathcal{H}_{{\bf
n},k}(x)d\sigma_k(x)}{Q_{{\bf n},k-1}(x) Q_{{\bf
n},k+1}(x)}=0.
\]
Separating this integral in two and using \eqref{rel1} we obtain \eqref{rel4}. \hfill $\Box$

\medskip

For ${\bf n} \in {\mathbb{Z}}_+^{m}(\bullet)$ and $\ell \in \{1,\ldots,m\}$ we define ${\bf n}^\ell$ as was done above. Though ${\bf n}^\ell$ need not belong to ${\mathbb{Z}}_+^{m}(\bullet)$, we can define the corresponding functions $\Psi_{{\bf n}^\ell,k}$.

\begin{teo} \label{teo:7} Assume that ${\bf n} \in {\mathbb{Z}}_+^{m}(\bullet)$ and $\ell \in \{1,\ldots,m\}$. For each $k=0,\ldots,m-1$,   ${\Psi}_{{\bf n}^\ell,k}$ has at most   $n_{k+1}+\cdots+n_m+1$ zeros in  ${\mathbb{C}} \setminus \Delta_{k}$ and at least $n_{k+1}+\cdots+n_m$ sign changes in $\Delta_{k+1}$. Therefore, all its zeros are real and simple. The zeros of ${\Psi}_{{\bf n},k}$ and ${\Psi}_{{\bf n}^\ell,k}$ in ${\mathbb{C}} \setminus \Delta_{k}$ interlace.
\end{teo}

{\bf Proof.} The proof is similar to that of Theorem \ref{teo:4} so we will not dwell into details. For each $k=0,\ldots,m-1$ and $A,B \in \mathbb{R}, |A| + |B| > 0$ define
\[ G_{{\bf n},k} = A \Psi_{{\bf n},k} + B \Psi_{{\bf n}^\ell,k}.
\]
Then
\[
\int G_{{\bf n},k-1}(x)\left(p_k(x) + \sum_{j=k+1}^m p_j(x) \widehat{s}_{k+1,j}(x)\right) d \sigma_k(x)  =0,\qquad
\deg p_j \leq n_j-1 \qquad j=k,\ldots,m.
\]
From here, it follows that there exists a monic polynomial $W_{{\bf n},k+1}, \deg W_{{\bf n},k+1} \geq N_{k+1} = n_{k+1} + \cdots, n_m, N_{m+1} = 0$,
whose roots are the zeros of $G_{{\bf n},k}$ in $\mathbb{C} \setminus \Delta_k$ such that
\[ \int  x^\nu  {G_{{\bf n},k-1}(x)} \frac{d\sigma_k(x)}{W_{{\bf n},k+1}(x)} = 0, \qquad \nu = 0,\ldots,n_k + \deg W_{{\bf n},k+1} -1.
\]
These relations imply that ${G}_{{\bf n},k}$ has at most   $n_{k+1}+\cdots+n_m+1$ zeros in  ${\mathbb{C}} \setminus \Delta_{k}$ and at least $n_{k+1}+\cdots+n_m$ sign changes in $\Delta_{k+1}$. In particular, this is true for $\Psi_{{\bf n}^\ell,k}$. So, all the zeros of ${G}_{{\bf n},k}$
(and $\Psi_{{\bf n}^\ell,k}$) in  ${\mathbb{C}} \setminus \Delta_{k}$ are real and simple. The interlacing is proved following the same arguments as in Theorem \ref{teo:4} The details are left to the reader. \hfill $\Box$

\section{Weak asymptotic.}

Following standard techniques,  the weak asymptotic for type I and type II Hermite-Pad\'e polynomials is derived using arguments from potential theory. We will briefly summarize what is needed.

\subsection{Preliminaries from potential theory.}

Let $E_k,\, k=1\ldots,m,$ be (not necessarily distinct)
compact subsets of the real line and
\[ \mathcal{C} = (c_{j,k}), \qquad
1\leq j,k \leq m,
\] a real, positive definite, symmetric
matrix of order $m$. $\mathcal{C}$ will be called the
interaction matrix. Let
$\mathcal{M}_1(E_k) $ be the subclass of probability measures in
$\mathcal{M}(E_k).$ Set
\[\mathcal{M}_1= \mathcal{M}_1(E_{1})
\times \cdots \times \mathcal{M}_1(E_{m})  \,.
\]

Given a vector measure $\vec{\mu}=(\mu_{1},\ldots,\,\mu_{m}) \in
\mathcal{M}_1$ and $j= 1,\ldots,m,$ we define the combined
potential
\[W^{\vec{\mu}}_j(x) = \sum_{k=1}^{m} c_{j,k}
V^{\mu_k}(x) \,,
\]
where
\[ V^{\mu_k}(x) := \int \log \frac{1}{|x-t|} \,d\mu_k(t)\,,
\]
denotes the standard logarithmic potential of $\mu_k$. Set
\[ \omega_j^{\vec{\mu}} := \inf \{W_j^{\vec{\mu}}(x): x \in E_j\} \,, \quad
j=1,\ldots,m\,.
\]

It is said that $\sigma \in \mathcal{M}(\Delta)$ is regular, and we write $\sigma \in \mbox{\bf Reg}$, if
\[ \lim \gamma_n^{1/n} = \frac{1}{\mbox{cap}(\mbox{supp}(\sigma))},
\]
where $\mbox{cap}(\mbox{supp}(\sigma) )$ denotes the logarithmic capacity of $\mbox{supp}(\sigma)$ and $\gamma_n$ is the leading coefficient of the (standard) $n$-th orthonormal polynomial with respect to $\sigma$. See \cite[Theorems 3.1.1, 3.2.1]{stto} for different equivalent forms of defining regular measures and its basic properties. In connection with regular measures it is frequently convenient that the support of the measure be regular. A compact set $E$ is said to be regular when the Green's function, corresponding to the unbounded connected component of $\mathbb{C} \setminus E$, with singularity at $\infty$ can be extended continuously to $E$.

In Chapter 5 of \cite{NiSo}  the authors prove (we state the
result in a form convenient for our purpose).

\begin{lemma} \label{niksor} Assume that the compact sets
$E_k,k=1,\ldots,m,$ are regular. Let $\mathcal{C}$ be a real, positive definite, symmetric
matrix of order $m$. If there exists $\vec{\lambda}  =
(\lambda_{1},\ldots,\lambda_{m})\in
\mathcal{M}_1$ such that for each $j=1,\ldots,m$
\[
W_j^{\vec{\lambda}} (x) = \omega_j^{\vec{\lambda}}\,, \qquad x
\in \supp{\lambda_j}\,,
\]
then $\vec{\lambda}$ is unique. Moreover, if $c_{j,k} \geq 0$
when $E_j \cap E_k \neq \emptyset$, then $\vec{\lambda}$ exists.
\end{lemma}

\medskip

For details on how Lemma \ref{niksor} is derived from  \cite[Chapter
5]{NiSo} see  \cite[Section 4]{bel}. The vector measure $\vec{\lambda}$
is called the equilibrium solution for the vector potential problem determined by the
interaction matrix $\mathcal{C}$ on the system of compact sets
$E_j\,, j = 1,\ldots,m$ and $\omega^{\vec{\lambda}}:= (\omega_1^{\vec{\lambda}},\ldots,\omega_m^{\vec{\lambda}})$ is the vector equilibrium constant. There are other characterizations of the equilibrium measure and constant but we will not dwell into that because they will not be used and their formulation requires introducing additional notions and notation.

\medskip

We also need

\begin{lemma} \label{lemextremal}
Let $E \subset \mathbb{R}$ be a regular compact set and $\phi$ a continuous function
on $E$. Then, there exists a unique $\lambda \in
\mathcal{M}_1(E)$ and a constant $w$ such that
\[
V^{\lambda}(z)+\phi(z) \left\{ \begin{array}{l} \leq
w,\quad z
\in \supp{\lambda} \,, \\
\geq w, \quad z \in  E\,.
\end{array} \right.
\]
\end{lemma}

In particular, equality takes place on all $\supp{\lambda}$.
If the compact set $E$ is not regular with respect to the
Dirichlet problem, the second part of the statement is true
except on a set $e$ such that $\mbox{cap}(e) =0.$ Theorem I.1.3 in
\cite{ST} contains a proof of this lemma in this context. When $E$
is regular, it is well known that this inequality except on a set
of capacity zero implies the inequality for all points in the set
(cf. Theorem I.4.8 from \cite{ST}). $\lambda$ is called the
equilibrium measure in the presence
of the external field $\phi$ on $E$ and $w$ is  the equilibrium
constant.

\medskip
As usual, a sequence of measures $(\mu_n)$ supported on a compact set $E$ is said to converge to a measure $\mu$ in the weak star topology if for every continuous function $f$ on $E$ we have
\[ \lim_{n} \int f d\mu_n = \int f d \mu.
\]
We write $*\lim_n \mu_n = \mu$. Given a polynomial $Q$ of degree $n$ we denote
\[ \mu_{Q} = \frac{1}{n} \sum_{Q(x) = 0} \delta_x
\]
where $\delta_x$ is the Dirac measure with mass $1$ at point $x$. In the previous sum, each zero of $Q$ is repeated taking account of its multiplicity. The measre $\mu_Q$ is usually called the normalized zero counting measure of $Q$.

One last ingredient needed is a result which relates the asymptotic zero distribution of polynomials orthogonal with respect to varying measures with the solution of a vector equilibrium problem in the presence of an external field contained in Lemma \ref{lemextremal}.  Different
versions of it appear in \cite{gora},  and
\cite{stto}. In \cite{gora}, it was proved assuming that
$\supp{\sigma}$ is an interval on which $\sigma'
> 0$ a.e.  Theorem 3.3.3 in
\cite{stto}  does not cover the type of
external field we need to consider. As stated here, the proof appears in \cite[Lemma 4.2]{FLLS}.

\begin{lemma}\label{gonchar-rakhmanov}
Assume that $\sigma \in \mbox{\bf Reg}$ and $\supp{\sigma} \subset
\mathbb{R}$  is regular.   Let $\{\phi_n\}, n \in \Lambda \subset
\mathbb{Z}_+,$ be a sequence of positive continuous functions on
$\supp{\sigma}$ such that
\begin{equation} \label{eq:phi}
\lim_{n\in \Lambda}\frac{1}{2n}\log\frac{1}{|\phi_n(x)|}= \phi(x)
> -\infty ,
\end{equation}
uniformly on $\supp{\sigma}$. Let  $(q_n), n \in \Lambda,$ be
a sequence of monic polynomials such that $\deg q_n = n$ and
\[
\int x^k q_n(x)\phi_n(x)d\sigma(x)=0,\qquad k=0,\ldots, n-1.
\]
Then
\begin{equation} \label{eq:18}
*\lim_{n \in \Lambda}\mu_{q_n} = \lambda,
\end{equation}
and
\begin{equation} \label{eq:19}
\lim_{n\in \Lambda}\left(\int |q_n(x)|^2\phi_n(x)
d\sigma(x)\right)^{1/{2n}}= e^{-w},
\end{equation}
where $\lambda$ and $w$ are the  equilibrium measure and
equilibrium constant in the presence of the external field $\phi$
on $\supp{\sigma}$ given by Lemma \ref{lemextremal}. We also have
\begin{equation} \label{eq:H}
\lim_{n \in \Lambda} \left(\frac{|q_n(z)|}{\|q_n
\phi_n^{1/2}\|_E}\right)^{1/n} = \exp{(w -
V^{\lambda}(z))}, \qquad  \mathcal{K} \subset  \mathbb{C} \setminus \Delta,
 \end{equation}
 where $\|\cdot\|_E$ denotes the uniform norm on $E$ and $\Delta$ is the smallest interval containing $\supp{\sigma}$.
\end{lemma}

\subsection{Weak asymptotic behavior for type II}
In the proof of the asymptotic zero distribution of the polynomials $Q_{{\bf n},j}$ we take $E_j = \supp{\sigma_j}$.  We need to specify the sequence of multi-indices for which the result takes place and the relevant interaction matrix for the vector equilibrium problem which arises.

Let $\Lambda =
\Lambda(p_{1},\ldots,p_{n})  \subset
\mathbb{Z}_+^{m}(\bullet)$ be an infinite sequence of distinct multi-indices
such that
\begin{equation} \label{indicesdec}
\lim_{ {\bf n} \in \Lambda} \frac{n_{j}}{|{\bf n}|} = p_{j}
\in (0,1), \qquad j=1,\ldots,m.
\end{equation}
Obviously, $p_{1}\geq \cdots\ge p_{m} $ and $\sum_{j=1}^{m} p_{j}=  1$. Set
\[
P_j=\sum_{k=j}^{m} p_{k},\qquad j=1,\ldots,m.
\]

Let us define the interaction matrix $\mathcal{C}_\mathcal{N}$ which is
relevant in the next result. Set
\begin{equation}\label{matriz}
\mathcal{C}_{\mathcal{N}}:=
\begin{pmatrix}
P_{1}^2 & -\frac{P_{1}P_{2}}{2} & 0 &  \cdots  & 0\\
-\frac{P_{1}P_{2}}{2} & P_{2}^2 &
-\frac{P_{2}P_{3}}{2} & \cdots  & 0\\
0 & -\frac{P_{2}P_{3}}{2} & P_{3}^2
& \cdots  & 0 \\
\vdots &\vdots&\vdots&\ddots &\vdots\\
0 & 0& 0 &\cdots  & P_{m}^2
\end{pmatrix}.
\end{equation}
This matrix satisfies all the assumptions of Lemma \ref{niksor} on
the compact sets $E_j = \supp(\sigma
_j), j=1,\ldots,m,$ including $c_{j,k} \geq
0$ when $E_j \cap E_k \neq \emptyset$  and it is positive definite because the principal
section $(\mathcal{C}_\mathcal{N})_r, r=1,\ldots,m$ of $\mathcal{C}_\mathcal{N}$
satisfies
\begin{displaymath}
\det(\mathcal{C}_\mathcal{N})_r=P_{1}^2\cdots P_{r}^2 \det
\begin{pmatrix}
1 & -\frac{1}{2} & 0 & \cdots & 0 & 0 \\
-\frac{1}{2} & 1 & -\frac{1}{2}  & \cdots & 0 & 0 \\
0 & -\frac{1}{2} &
1  & \cdots & 0 & 0 \\
\vdots & \vdots & \vdots   & \ddots & \vdots & \vdots\\
0 & 0 & 0 & \cdots & 1 & -\frac{1}{2}\\
0 & 0 & 0  & \cdots & -\frac{1}{2} & 1
\end{pmatrix}_{r\times r} > 0.
\end{displaymath}
Let $\vec{\lambda}(\mathcal{C}_\mathcal{N}) = (\lambda_1,\ldots,\lambda_m)$ be the solution of the corresponding vector equilibrium problem stated in Lemma \ref{niksor}.

The next result, under more restrictive conditions on the measures but in the framework of so called Nikishin systems on a graph tree is contained in \cite{GRS}. We have practically reproduced their arguments which incidentally can also be adapted to the study of so called mixed type Hermite-Pad\'e approximation in which the definition contains a mixture of type I and type II interpolation conditions. For details see \cite{FLLS}.

 \begin{teo} \label{teo4} Let $\Lambda$ be a sequence of multi-indices verifying \eqref{indicesdec}. Assume that $\sigma_j \in \mbox{\bf Reg}$ and $\supp{\sigma_j} = E_j$ is regular for each $j=1,\ldots,m$. Then,
\begin{equation} \label{weak}
*\lim_{{\bf n}\in \Lambda} \mu_{Q_{{\bf n},j}} =  {\lambda}_{j}, \qquad j=1,\ldots,m.
\end{equation}
where $\vec{\lambda} = (\lambda_1,\ldots,\lambda_m) \in
\mathcal{M}_1$ is the vector equilibrium measure determined by the
matrix $\mathcal{C}_{\mathcal{N}}$ on the system of compact
sets $E_j, j=1,\ldots,m$.  Moreover,
\begin{equation} \label{eq:4*}
\lim_{{\bf n}\in \Lambda} \left|\int  {Q_{{\bf n},j}^2(x)}\frac{\mathcal{H}_{{\bf n},j}(x)\,{\rm d}\sigma_j(x)}{Q_{{\bf n},j-1}(x)Q_{{\bf n},j+1}(x)} \right|^{1/2|{\bf n}|} = \exp\left(-\sum_{k=1}^{j}
{\omega_k^{\vec{\lambda}}/P_k} \right) \,,
\end{equation}
where  $\omega^{\vec{\lambda}} = (\omega_1^{\vec{\lambda}},\ldots,\omega_m^{\vec{\lambda}})$ is the vector equilibrium constant. For $j=1,\ldots,m$
\begin{equation} \label{asint1}
\lim_{{\bf n}\in \Lambda} |\Psi_{{\bf n},j}(z)|^{1/|{\bf n}|} = \exp\left(P_j V^{{\lambda}_{j}}(z)- P_{j+1}V^{ {\lambda}_{j+1}}(z)-2 \sum_{k=1}^{j}
{\omega_k^{\vec{\lambda}}/P_k} \right)
\end{equation}
uniformly on compact subsets of $\mathbb{C} \setminus (\Delta_j \cup \Delta_{j+1})$ where $\Delta_{m+1} = \emptyset$ and the term with $P_{m+1}$ is dropped when $j=m$.
\end{teo}
{\bf Proof.} The unit ball in the cone of positive Borel measures
is weak star compact; therefore, it is sufficient to show that
each one of the sequences of measures $(\mu_{Q_{{\bf n},j}})$,
${\bf n}\in\Lambda$, $j=1,\ldots,m,$ has only one
accumulation point which coincides with the corresponding
component of the vector equilibrium measure $\vec{\lambda}$ determined by the
matrix $\mathcal{C}_{\mathcal{N}}$ on the system of compact
sets $E_j, j=1,\ldots,m$.

\medskip

Let
$\Lambda' \subset \Lambda$ be  such
that for each $j=1,\ldots,m$
\[
*\lim_{{\bf n}\in\Lambda'}\mu_{Q_{{\bf n},j}}=\mu_j.
\]
Notice that $\mu_j\in\mathcal{M}_1(E_j)$, $j=1,\ldots,m$. Taking into account that all the zeros of $Q_{{\bf n},j}$ lie in $\Delta_j$, it follows that
\begin{equation}\label{conv-Qnj}
\lim_{{\bf n}\in\Lambda'}|Q_{{\bf n},j}(z)|^{1/|{\bf n}|}=\exp(- P_jV^{\mu_j}(z)),
\end{equation}
uniformly on compact subsets of $\mathbb{C} \setminus\Delta_j$.

\medskip

When $k=1$, \eqref{rel3} reduces to
\[ \int x^{\nu} Q_{{\bf n},1}(x) \frac{{\rm d} \sigma_{1}(x)}{|Q_{{\bf
n},2}(x)|} = 0\,, \qquad \nu=0,\ldots,|{\bf n}|-1.
\]
According to (\ref{conv-Qnj})
\[ \lim_{{\bf n} \in \Lambda'} \frac{1}{2|{\bf n}|}\log|Q_{{\bf n},2}(x)| =
-\frac{P_2}{2 } V^{\mu_{2}}(x)\,,
\]
uniformly on $\Delta_{2}$. Using Lemma \ref{gonchar-rakhmanov},
it follows that $\mu_{1}$ is the unique solution of the extremal
problem
\begin{equation} \label{eq:1}
V^{\mu_{1}}(x) - \frac{P_2}{2}V^{\mu_{2}}(x)
\left\{
\begin{array}{l} = \omega_{1},\quad x
\in \supp {\mu_{1}} \,, \\
\geq \omega_{1}, \quad x \in E_{1} \,,
\end{array} \right.
\end{equation}
and
\begin{equation} \label{eq:2}
\lim_{{\bf n} \in \Lambda'} \left|\int \frac{Q_{{\bf n},1}^2(x)}{|Q_{{\bf n},2}(x)|}{\rm d}\sigma_{1}(x)\right|^{1/2|{\bf n}|} =
e^{-\omega_{1}}\,.
\end{equation}

Using induction on increasing values of $j$, let us show that for all
$j = 1,\ldots,m$
\begin{equation} \label{eq:3*}
V^{\mu_j}(x) - \frac{P_{j-1}}{2 P_j}V^{\mu_{j-1}}(x) -
\frac{P_{j+1}}{2 P_j}V^{\mu_{j+1}}(x) +
 \frac{P_{j-1}}{P_j}\omega_{j-1} \left\{
\begin{array}{l} = \omega_j,\quad x
\in \supp{\mu_j} \,, \\
\geq \omega_j, \quad x \in E_j \,,
\end{array} \right.,
\end{equation}
(when $j=1$ or $j=m$ the terms with $P_0$ and $P_{m+1}$ do not appear,)    and
\begin{equation} \label{eq:4**}
\lim_{n \in \Lambda'} \left|\int  {Q_{{\bf n},j}^2(x)} \frac{|\mathcal{H}_{{\bf n},j}(x)| {\rm d}\sigma_j(x)}{|Q_{{\bf n},j-1}(x)Q_{{\bf n},j+1}(x)|} \right|^{1/2N_{{\bf n},j}} = e^{-\omega_j}\,,
\end{equation}
where $Q_{n,0} \equiv Q_{n,m+1} \equiv 1$ and $N_{{\bf n},j} = n_j + \cdots + n_m$. For $j = 1$ these relations
are non other than (\ref{eq:1})-(\ref{eq:2}) and the initial
induction step is settled. Let us assume that the statement is
true for $j -1\in \{1,\ldots,m-1\}$ and let us prove it for
$j$.

\medskip

Taking acount of the fact that $Q_{{\bf n},j-1},Q_{{\bf n},j+1}$  and $\mathcal{H}_{{\bf n},j}$ have constant sign on $\Delta_j$, for $j=1,\ldots, m$, the orthogonality relations  \eqref{rel3} can be expressed as
\[ \int x^{\nu} Q_{{\bf n},j}(x)  \frac{|\mathcal{H}_{{\bf n},j}(x)| {\rm d}\sigma_j(x)}{|Q_{{\bf n},j-1}(x)Q_{{\bf
n},j+1}(x)|} = 0\,, \qquad \nu=0,\ldots,N_{{\bf n},j}-1\,,
\]
and using \eqref{rel4} it follows that
\[ \int x^{\nu} Q_{{\bf n},j}(x) \left| \int \frac{Q_{{\bf n},j-1}^2(t)}{|x-t|} \frac{|\mathcal{H}_{{\bf n},j-1}(t)|{\rm d}\sigma_{j-1}(t)}{ |Q_{{\bf n},j-2}(t)Q_{{\bf n},j}(t)| }\right|\frac{{\rm d}\sigma_j(x)}{|Q_{{\bf n},j-1}(x)Q_{{\bf n},j+1}(x)|} = 0\,,
\]
for $\nu=0,\ldots,N_{{\bf n},j}-1\,.$

\medskip

Relation (\ref{conv-Qnj}) implies that
\begin{equation}\label{eq:5}
\lim_{n \in \Lambda'} \frac{1}{2N_{{\bf n},j}}\log|Q_{{\bf n},j-1}(x)Q_{{\bf n},j+1}(x)| = -
\frac{P_{j-1}}{2 P_j}V^{\mu_{j-1}}(x) -
\frac{P_{j+1}}{2 P_j}V^{\mu_{j+1}}(x)\,,
\end{equation}
uniformly on $\Delta_j.$ (Since $Q_{n,0}\equiv 1$, when
$j=1$ we only get the second term on the right hand side of
this limit.)

\medskip

Set
\begin{equation} \label{Knj} K_{n,j-1} := \left|\int  {Q_{{\bf n},j-1}^2(t)}  \frac{ |{\mathcal{H}}_{{\bf n},j-1}(t)|{\rm d}\sigma_{j-1}(t)}{|Q_{{\bf n},j-2}(t)Q_{{\bf n},j}(t)|}
 \right|^{-1/2}.
\end{equation}
It follows that for $x \in \Delta_j$
\[ \frac{1}{\delta_{j-1}^*K_{{\bf n},j-1}^2} \leq \left|\int \frac{Q_{{\bf n},j-1}^2(t)}{|x-t|}\frac{ |{\mathcal{H}}_{{\bf n},j-1}(t)|{\rm d}\sigma_{j-1}(t)}{|Q_{{\bf n},j-2}(t)Q_{{\bf n},j}(t)|}
\right|\leq \frac{1}{\delta_{j-1}K_{n,j-1}^2},
\]
where $0 < \delta_{j-1} = \min\{|x-t|: t \in \Delta_{j-1}, x \in
\Delta_j\} \leq \max\{|x-t|: t \in \Delta_{j-1}, x \in \Delta_j\}
= \delta_{j-1}^* < \infty.$ Taking into consideration these
inequalities, from the induction hypothesis, we obtain that
\begin{equation} \label{eq:6}
\lim_{n \in \Lambda'} \left|\int \frac{Q_{{\bf n},j-1}^2(t)}{|x-t|}\frac{ |{\mathcal{H}}_{{\bf n},j-1}(t)|{\rm d}\sigma_{j-1}(t)}{|Q_{{\bf n},j-2}(t)Q_{{\bf n},j}(t)|}\right|^{1/2N_{{\bf n},j}} = e^{- P_{j-1}\omega_{j-1}/P_j }.
\end{equation}

\medskip

Taking (\ref{eq:5}) and (\ref{eq:6}) into account, Lemma
\ref{gonchar-rakhmanov} yields that $\mu_j$ is the unique solution
of the extremal problem (\ref{eq:3*}) and
\[  \lim_{n \in \Lambda'} \left|\int   \int \frac{Q_{{\bf n},j-1}^2(t)}{|x-t|}\frac{ |{\mathcal{H}}_{{\bf n},j-1}(t)|{\rm d}\sigma_{j-1}(t)}{|Q_{{\bf n},j-2}(t)Q_{{\bf n},j}(t)|}
\frac{Q_{{\bf n},j}^2(x){\rm d}\sigma_j(x)}{|Q_{{\bf n},j-1}(x)Q_{{\bf n},j+1}(x)|}\right|^{1/2N_{{\bf n},j} } = e^{-\omega_j}.
\]
According to \eqref{rel4} the previous formula reduces to \eqref{eq:4**}. We have concluded the induction.

\medskip

Now, we can rewrite (\ref{eq:3*}) as
\begin{equation} \label{eq:a1}
 P_j^2V^{\mu_j}(x) - \frac{P_jP_{j-1}}{2}V^{\mu_{j-1}}(x) -
\frac{P_jP_{j+1}}{2 }V^{\mu_{j+1}}(x)  \left\{
\begin{array}{l} = \omega_j',\quad x
\in \supp{\mu_j} \,, \\
\geq \omega_j', \quad x \in E_j \,,
\end{array} \right.
\end{equation}
for $j=1,\ldots,m$, where
\begin{equation}\label{eq:c}
\omega_j' =  P_j^2\omega_j -  P_jP_{j-1}\omega_{j-1}, \qquad (\omega_{0} = 0).
\end{equation}
(Recall that the terms with $V^{\mu_0}$ and $V^{\mu_{m+1}}$ do not appear when $j=0$ and $j=m$, respectively.)
By Lemma \ref{niksor},
$\vec{\lambda} = (\mu_{1},\ldots,\mu_{m})$ is the solution of the equilibrium  problem determined by the
interaction matrix $\mathcal{C}_{\mathcal{N}}$ on the system of compact sets
$E_j\,, j = 1,\ldots,m$ and $\omega^{\vec{\lambda}} = (\omega_{1}',\ldots,\omega_{m}') $  is the corresponding vector equilibrium constant.  This is
for any convergent subsequence; since the equilibrium problem does not depend on the sequence of indices $\Lambda'$  and the solution is unique we  obtain the limits in \eqref{weak}. By the same token, the limit in \eqref{eq:4**} holds true over the whole sequence of indices $\Lambda$. Therefore,
\begin{equation} \label{eq:4***}
\lim_{n \in \Lambda} \left|\int  {Q_{{\bf n},j}^2(x)} \frac{|\mathcal{H}_{{\bf n},j}(x)| {\rm d}\sigma_j(x)}{|Q_{{\bf n},j-1}(x)Q_{{\bf n},j+1}(x)|} \right|^{1/2|{\bf n}|} = e^{-P_j\omega_j}\,.
\end{equation}

\medskip

From (\ref{eq:c}) it follows that
$\omega_{1} = \omega_{1}^{\vec{\lambda}}$ when
$j=1.$  Suppose that $ P_{j-1}\omega_{j-1} = \sum_{k=1}^{j-1}
 {\omega_{k}^{\vec{\lambda}}}/P_k$ where $j-1 \in
\{1,\ldots,m-1\}$. Then, according to (\ref{eq:c})
\[ P_j \omega_j = {\omega_j^{\vec{\lambda}}}/P_j  +
 P_{j-1}\omega_{j-1} = \sum_{k=1}^{j}
 {\omega_{k}^{\vec{\lambda}}}/P_k
\]
and (\ref{eq:4*}) immediately follows using \eqref{eq:4***}.

For $j \in \{1,\ldots,m\}$, from \eqref{rel4}, we have
\begin{equation} \label{eq:7} \Psi_{{\bf n},j}(z) = \frac{Q_{{\bf n},j+1}(z)}{Q_{{\bf n},j}(z)}\int \frac{ Q^2_{{\bf n},j}(x)}{z-x} \frac{\mathcal{H}_{{\bf n},j}(x){\rm d}\sigma_{j}(x)}{Q_{{\bf n},j-1}(x)Q_{{\bf n},j+1}(x)},
\end{equation}
where $Q_{{\bf n},0}\equiv Q_{{\bf n},m+1} \equiv 1$. Now, \eqref{weak} implies
\begin{equation} \label{eq:7*} \lim_{{\bf n} \in \Lambda}\left|\frac{Q_{{\bf n},j+1}(z)}{Q_{{\bf n},j}(z)}\right|^{ {1}/|{\bf n}|}=\exp\left(
P_j V^{{\lambda}_{j}}(z)- P_{j+1}V^{ {\lambda}_{j+1}}(z)\right),
\end{equation}
uniformly on compact subsets of $\mathbb{C} \setminus (\Delta_j \cup
\Delta_{j+1})$ (we also use that the zeros of $Q_{n,j}$ and $Q_{n,j+1}$ lie in $\Delta_j$ and $\Delta_{j+1}$, respectively). It
remains to find the $|{\bf n}|$-th root asymptotic behavior of
the integral.

Fix a compact set $\mathcal{K} \subset \mathbb{C} \setminus
\Delta_{j}.$
It is not difficult to prove that (for the definition of
$K_{n,j}$ see \eqref{Knj})
\[ \frac{C_1}{K_{{\bf n},j}^2} \leq \left|\int \frac{ Q^2_{{\bf n},j}(x)}{z-x} \frac{\mathcal{H}_{{\bf n},j}(x){\rm d}\sigma_{j}(x)}{Q_{{\bf n},j-1}(x)Q_{{\bf n},j+1}(x)}\right| \leq \frac{C_2}{K_{{\bf
n},j}^2} \,,
\]
where
\[ C_1 = \frac{\min \{ \max\{|u-x|,|v|: z = u+iv\}: z \in \mathcal{K}, x
\in \Delta_{j}\}}{ \max\{|z-x|^2: z \in \mathcal{K}, x \in
\Delta_{j}\}} > 0
\]
and
\[ C_2 = \frac{1}{\min\{|z-x|: z \in \mathcal{K}, x \in
\Delta_{j}\}} < \infty.
\]
Taking into account (\ref{eq:4*})
\begin{equation}\label{eq:9}
\lim_{{\bf n} \in \Lambda} \left|\int \frac{ Q^2_{{\bf n},j}(x)}{z-x} \frac{\mathcal{H}_{{\bf n},j}(x){\rm d}\sigma_{j}(x)}{Q_{{\bf n},j-1}(x)Q_{{\bf n},j+1}(x)}\right|^{1/|{\bf n}|} = \exp\left(-2 \sum_{k=1}^{j}
{\omega_k^{\vec{\lambda}}/P_k} \right) \,.
\end{equation}
From \eqref{eq:7},  \eqref{eq:7*}, and (\ref{eq:9}), we obtain \eqref{asint1} and we
are done.
\hfill $\Box$

\begin{rem} In the case of type I Hermite Pad\'e approximation, asymptotic formulas for the forms $\mathcal{A}_{{\bf n},j}$ and the polynomials $A_{{\bf n},j}$  can be obtained following  arguments similar to those employed above, see \cite{nik2} or \cite{NiSo}. Basically, all one has to do is replace the use of the formulas in Proposition \ref{pro1} by the ones in Proposition \ref{teo:5}. I recommend doing it as an exercise.
\end{rem}

\subsection{Application to Hermite-Pad\'e approximation}

The convergence of type II Hermite-Pad\'e approximants for the case of $m$ generating measures and interpolation conditions equally distributed between the different functions was obtained in \cite{Bus}. When $m=2$ the result was proved in \cite{Nik}. Other results for more general sequences of multi-indices and so called multipoint Hermite-Pad\'e approximation were considered in \cite{LF2, LF3}. For type I, the convergence was proved recently in \cite{LS}. Here we wish to show how the weak  asymptotic of the Hermite-Pad\'e polynomials allows to estimate the rate of convergence of the approximants. We restrict  to type II. The presentation follows closely the original result given in \cite{GRS}.

Consider the functions
\begin{equation}\label{remain} \Phi_{{\bf n},j}(z) := (Q_{\bf n} \widehat{s}_{1,j} - P_{{\bf n},j})(z) = \mathcal{O}(1/z^{n_j +1}), \qquad z \to \infty, \qquad j=1,\ldots,m.
\end{equation}
which are the remainders of the interpolation conditions defining the type II Hermite-Pad\'e approximants with respect to the multi-index $\bf n$ of the Nikishin system of functions $(\widehat{s}_{1,1},\ldots,\widehat{s}_{1,m})$.
Because of \eqref{remain}, $P_{{\bf n},j}$ is the polynomial part of the Laurent expansion at $\infty$ of $Q_{\bf n} \widehat{s}_{1,j}$. It is easy to check that
\[(Q_{\bf n} \widehat{s}_{1,j} - P_{{\bf n},j})(z) =  \int \frac{Q_{\bf n}(x) ds_{1,j}(x)}{z-x}, \qquad  P_{{\bf n},j}(z) = \int \frac{Q_{\bf n}(z) - Q_{\bf n}(x)}{z-x} d s_{1,j}(x).
\]
For example, this follows using Hermite's integral representation of $Q_{\bf n} \widehat{s}_{1,j} - P_{{\bf n},j}$, Cauchy's integral formula, and the Fubini theorem. According to the way in which $s_{1,j}$ is defined, we have
\[ \Phi_{{\bf n},j}(z) = \int \cdots \int  \frac{Q_{\bf n}(x_1) d\sigma_1(x_1)d\sigma_2(x_2)\cdots d\sigma_j(x_j)}{(z- x_1)(x_1-x_2)\cdots(x_{j-1} - x_{j})}.
\]
Notice that $\Phi_{{\bf n},1} = \Psi_{{\bf n},1}$. We wish to establish a connection between the functions $\Phi_{{\bf n},j}$ and $\Psi_{{\bf n},k}, 1\leq j,k \leq m$.

First let us present an interesting formula which connects $(s_{1,1},\ldots,s_{1,m}) = \mathcal{N}(\sigma_1,\ldots,\sigma_m)$ and  $(s_{m,m},\ldots,s_{m,1}) = \mathcal{N}(\sigma_m,\ldots,\sigma_1)$. When $1\leq k < j \leq m$, we denote
\[ s_{j,k} = \langle \sigma_j,\sigma_{j-1},\ldots,\sigma_k\rangle.
\]

\begin{lemma} \label{milagro}
For each $j=2,\ldots,m$,
\begin{equation} \label{miracle}
(\widehat{s}_{1,j} - \widehat{s}_{1,j-1}\widehat{s}_{j,j} + \widehat{s}_{1,j-2}\widehat{s}_{j,j-1} + \cdots + (-1)^{j-1} \widehat{s}_{1,1}\widehat{s}_{j,2} +(-1)^j \widehat{s}_{j,1})(z) \equiv 0,
\end{equation}
for all $z \in \mathbb{C} \setminus (\Delta_1 \cup \Delta_j)$.
\end{lemma}
{\bf Proof.} Notice that
\[ \widehat{s}_{1,j}(z) + (-1)^j \widehat{s}_{j,1}(z) = \int \cdots \int \frac{(x_1 - x_j)d\sigma_1(x_1)d\sigma_2(x_2)\cdots d\sigma_j(x_j)}{(z-x_1)(x_1-x_2) \cdots (x_{j-1}-x_j)(z-x_j)}.
\]
On the right hand side, use that $x_1-x_j = (x_1 - x_2) + (x_2 - x_3) + \cdots + (x_{j-1} - x_j)$ to separate the integral in a sum. In each one of the resulting integrals, the numerator cancels one of the factors in the denominator, and the integral splits in the product of two which easily identify  with the remaining terms in the formula. \hfill $\Box$

Now we can prove the connection formulas.

\begin{lemma} \label{connection}
We have that $\Psi_{{\bf n},1} = \Phi_{{\bf n},1}$ and for $j = 2,\ldots,m$
\begin{equation} \label{con1}
\Psi_{{\bf n},j}(z)=
\sum_{k=2}^{j} (-1)^{k} \widehat{s}_{j,k}(z)
 {\Phi}_{ {\bf n},k-1}(z) + (-1)^{j+1} {\Phi}_{ {\bf n},j}(z), \qquad z \in {\bbc}
\setminus  (\Delta_1 \cup \Delta_j) \,,
\end{equation}
and
\begin{equation} \label{con2}
\Phi_{{\bf n},j}(z)=
\sum_{k=2}^{j} (-1)^{k} \widehat{s}_{k,j}(z)
 {\Psi}_{ {\bf n},k-1}(z) + (-1)^{j+1} {\Psi}_{ {\bf n},j}(z), \qquad z \in {\bbc}
\setminus   (\cup_{k=1}^j \Delta_k) \,.
\end{equation}
\end{lemma}

{\bf Proof.} Obviously, $\Psi_{{\bf n},1} = \Phi_{{\bf n},1}$. Notice that formula \eqref{miracle} remains valid if the measures $\sigma_1,\ldots,\sigma_m$ are signed and finite. All what is needed  is that they are supported on intervals which are consecutively non intersecting. It is easy to see that
\[ \Phi_{{\bf n},j}(z) = \langle Q_{\bf n}\sigma_1,\sigma_2,\ldots,\sigma_j \widehat{\rangle} (z), \qquad j=2,\ldots,m.
\]
The symbol $\langle \cdot \widehat{\rangle}$ means taking the Cauchy transform of $\langle \cdot {\rangle}$. On the other hand
\[ \Psi_{{\bf n},j}(z) = \langle \sigma_j, \ldots,\sigma_2, Q_{\bf n}\sigma_1\widehat{\rangle}(z), \qquad j=2,\ldots,m.
\]
Taking into consideration the previous remarks, using formula \eqref{miracle} with $d\sigma_1$ replaced with $Q_{\bf n} d \sigma_1$, after trivial transformations we obtain \eqref{con1}.

The collection of formulas \eqref{con1} for $j=2,\ldots,m$ together with $\Phi_{{\bf n},1} = \Psi_{{\bf n},1}$ can be expressed in matrix form as follows
\[ (\Psi_{{\bf n},1},\ldots,\Psi_{{\bf n},m})^t = D (\Phi_{{\bf n},1},\ldots,\Phi_{{\bf n},m})^t,
\]
where $(\cdot)^t$ is the transpose of the vector $(\cdot)$ and $D$ is the $m\times m$ lower triangular matrix given by
\[ D :=
\left(
\begin{array}{ccccc}
1 & 0 & 0 & \cdots & 0 \\
\widehat{s}_{2,2} & -1 & 0 & \cdots & 0 \\
\widehat{s}_{3,2} & - \widehat{s}_{3,3} & 1 & \cdots & 0 \\
\vdots & \vdots & \vdots & \ddots & \vdots \\
\widehat{s}_{m,2} & - \widehat{s}_{m,3} & \widehat{s}_{m,4} & \cdots & (-1)^{m-1}
\end{array}
\right).
\]
Obviously $D$ is invertible and
\begin{equation} \label{inversion} (\Phi_{{\bf n},1},\ldots,\Phi_{{\bf n},m})^t = D^{-1} (\Psi_{{\bf n},1},\ldots,\Psi_{{\bf n},m})^t
\end{equation}
is the matrix form of the relations which express each function $\Phi_{{\bf n},j}$ in terms of $\Psi_{{\bf n},k}, k=1,\ldots,j$.
The matrix $D^{-1}$ is also lower triangular and it may be proved, using \eqref{miracle} in several ways, that
\[ D^{-1} :=
\left(
\begin{array}{ccccc}
1 & 0 & 0 & \cdots & 0 \\
\widehat{s}_{2,2} & -1 & 0 & \cdots & 0 \\
\widehat{s}_{2,3} & - \widehat{s}_{3,3} & 1 & \cdots & 0 \\
\vdots & \vdots & \vdots & \ddots & \vdots \\
\widehat{s}_{2,m} & - \widehat{s}_{3,m} & \widehat{s}_{4,m} & \cdots & (-1)^{m-1}
\end{array}
\right).
\]
Using \eqref{inversion} and the expression of $D^{-1}$ we obtain \eqref{con2}. \hfill $\Box$

\medskip

Let $\vec{\lambda} = (\lambda_1,\ldots,\lambda_m) \in
\mathcal{M}_1$ be the vector equilibrium measure determined by the
matrix $\mathcal{C}_{\mathcal{N}}$ on the system of compact
sets $E_j = \supp(\sigma_j), j=1,\ldots,m$. In the sequel we assume that the hypothesis of Theorem \ref{teo4} hold. For each $j=1, \ldots, m$, set
\[
U_j^{\vec{\lambda}}=P_j V^{{\lambda}_{j}}(z)- P_{j+1}V^{ {\lambda}_{j+1}}(z)-2 \sum_{k=1}^{j}
{\omega_k^{\vec{\lambda}}/P_k}\,,
\]
($V^{\widehat{\mu}_{m+1}} \equiv 0 $). Notice that in a
neighborhood of $z=\infty$ we have
\[
U^{\vec{\lambda}}_j(z) =  {\mathcal{O}}\left(p_{j} \log
\frac{1}{|z|}\right).
\]
The potentials of the components of the equilibrium measure  define continuous functions on all $\mathbb{C}$ (see the equilibrium equations). Thus, the functions $U^{\vec{\lambda}}_j$ are defined and continuous on all $\mathbb{C}$.

Fix $j \in \{1, \ldots m\}$. For  $k=1, \ldots ,j$ define the
regions
\[
D^j_k=\{z \in   {\mathbb{C}}  :
U^{\vec{\lambda}}_k(z)
> U^{\vec{\lambda}}_i(z), i= 1,\ldots,j \}.
\]
Some $D^j_k$ could be empty. Denote
\[
\xi_j(z)= \max \{ U^{\vec{\lambda}}_k(z): k=1, \ldots j \}.
\]

\begin{cor} \label{cor3} Under the assumptions of Theorem \ref{teo4}, for each $j=1,\ldots,m$, we have
\begin{equation}
\label{impasintotica} \lim_{{\bf n} \in \Lambda}
\left|\widehat{s}_{1,j}(z) - \frac{P_{{\bf n},j}(z)}{Q_{\bf n}(z)} \right|^{1/|{\bf n}|} = \exp (V^{\lambda_1}
+ \xi_j)(z)\,, \qquad z \in (\cup_{k=1}^j D^j_k) \setminus (\cup_{k=1}^{j+1} \Delta_k ) \,,
\end{equation}
and
\begin{equation} \label{eq:6.4} \limsup_{{\bf n} \in \Lambda}
\left|\widehat{s}_{1,j}(z) - \frac{P_{{\bf n},j}(z)}{Q_{\bf n}(z)} \right|^{1/|{\bf n}|} \leq \exp (V^{\lambda_1}
+ \xi_j)(z) \,, \qquad z \in \mathbb{C} \setminus (\cup_{k=1}^{j+1} \Delta_k )\,,
\end{equation}
uniformly on compact subsets   of the indicated regions. Moreover, $(V^{\lambda_1}
+ \xi_j)(z) < 0, z \in \mathbb{C} \setminus \Delta_1$ which implies that the sequence $({P_{{\bf n},j}}/{Q_{\bf n}}), {\bf n} \in \Lambda,$ converges to $\widehat{s}_{1,j}$ with geometric rate in  $\mathbb{C} \setminus (\cup_{k=1}^{j+1} \Delta_k )$.
\end{cor}

{\bf Proof.} By \eqref{asint1} and \eqref{con2}
we have that the following asymptotic formula takes place (notice that the functions $\widehat{s}_{k,j}$ are different from zero in ${\mathbb{C}} \setminus (\cup_{k=1}^j \Delta_k)$),
\[
\lim_{n \in
\Lambda}|{\Phi}_{ {\bf n},j}(z)|^{1/{|\bf n|}}=
\exp U^{\vec{\lambda}}_k(z), \qquad z \in D^j_k \setminus  (\cup_{k=1}^{j+1} \Delta_k),
\]
uniformly on compact subsets of the specified region. Then
\[ \lim_{n \in \Lambda}
| {\Phi}_{{\bf n},j}(z)|^{{1}/{|\bf n|}} = \exp
\xi_j(z)\,, \qquad z \in (\cup_{k=1}^j D^j_k) \setminus (\cup_{k=1}^{j+1} \Delta_k )\,,
\]
and
\[ \limsup_{n \in \Lambda}
|{\Phi}_{{\bf n},j}(z)|^{ {1}/{|\bf n|}} \leq \exp
\xi_j(z) \,, \qquad z \in \mathbb{C} \setminus (\cup_{k=1}^{j+1} \Delta_k )\,,
\]
uniformly on compact subsets of the specified region.

Formulas \eqref{impasintotica} and \eqref{eq:6.4} follow directly from
\[ \widehat{s}_{1,j}(z) - \frac{P_{{\bf n},j}(z)}{Q_{\bf n}(z)} = \frac{\Phi_{{\bf n},j}(z)}{Q_{\bf n}(z)},
\]
the asymptotic formulas given for $\Phi_{{\bf n},j}$,  and \eqref{weak}.

When $j=1$, we have
\[ (V^{\lambda_1} + \xi_1)(z) = 2 V^{\lambda_1}(z) - P_2 V^{\lambda_2}(z) -2\omega_1^{\vec{\lambda}} = 2(W_1^{\vec{\lambda}}(z) - \omega_1^{\vec{\lambda}}).
\]
According to \eqref{eq:1}, $W_1^{\vec{\lambda}}(z) - \omega_1^{\vec{\lambda}} \equiv 0, x \in \supp(\lambda_1)$. On the other hand, $W_1^{\vec{\lambda}}(z) - \omega_1^{\vec{\lambda}} $ is subharmonic in ${\mathbb{C}} \setminus \supp(\lambda_1)$ and tends to $-\infty$ as $z \to \infty$. By the maximum principle for subharmonic functions $W_1^{\vec{\lambda}}(z) - \omega_1^{\vec{\lambda}} < 0, z \in {\mathbb{C}} \setminus \supp(\lambda_1)$ (equality cannot occur at any point of this region because it would imply that $W_1^{\vec{\lambda}}(z) - \omega_1^{\vec{\lambda}} \equiv  0$ which is impossible).

Let us assume that $(V^{\lambda_1}
+ \xi_{j-1})(z) < 0, z \in \mathbb{C} \setminus \Delta_1$, where $j \in
\{2,\ldots,m\}$ and let us prove that $(V^{\lambda_1}
+ \xi_{j})(z) < 0, z \in \mathbb{C} \setminus \Delta_1$.
Obviously,
\[ \xi_j(z) = \max\{\xi_{j-1}(z),U_j(z)\}.\]
Consider the
difference
\[ U_j(z) - U_{j-1}(z) = 2(W_j^{\vec{\lambda}}(z) -
w_j^{\vec{\lambda}})/P_j =
O\left((p_{j}-p_{j-1})\log({1}/{|z|})\right), z \to
\infty.
\]

If $p_{j}=p_{j-1} = 0$ then $W_j^{\vec{\lambda}}(z) -
w_j^{\vec{\lambda}}$ is subharmonic in $\overline{\mathbb{C}}
\setminus \supp(\lambda_j)$ (at $\infty$ it is finite) and equals zero on
$\supp( \lambda_j)$. Hence, $U_j(z) \leq U_{j-1}(z) \leq
\xi_{j-1}(z)$ on $ {\mathbb{C}} \setminus
\supp(\lambda_j)$. Therefore, using the equilibrium
condition, $U_j(z) = U_{j-1}(z)$ on $\Delta_j$ and $U_j(z) <
U_{j-1}(z)$ on $ {\mathbb{C}} \setminus \Delta_j$. In this
case, $\xi_j(z) = \xi_{j-1}(z), z \in \mathbb{C} \setminus \Delta_1,$ and the conclusion
follows from the induction hypothesis.

If $p_{j} < p_{j -1}$, in a neighborhood of $\infty$ we
have $U_j(z) > U_{j-1}(z)$ since
$(p_{j}-p_{j-1})\log({1}/{|z|}) \to +\infty$ as $z
\to \infty.$ Let $\Gamma = \{z \in \mathbb{C}: U_j(z) = U_{j-1}(z)\}$. This
set contains $\supp( {\lambda}_j)$ and divides $\mathbb{C} \setminus \Delta_1$ in two
domains $\Omega_1 = \{z \in\mathbb{C} \setminus \Delta_1: U_j(z) > U_{j-1}(z)\}$, which
contains $z= \infty$, and $\Omega_2 = \{z \in \mathbb{C} \setminus \Delta_1: U_j(z) <
U_{j-1}(z)\}$. Since $U_{j-1}(z) \leq \xi_{j-1}(z)$,   on
$\Omega_2 \cup \Gamma$ we have that $\xi_{j-1}(z) = \xi_{j}(z)$
and thus  $(V^{\lambda_1} + \xi_{j}) < 0$.   On
$\Omega_1$ the function $V^{\lambda_1} + U_{j}$ is
subharmonic and on its boundary $\Gamma$ equals
$V^{\lambda_1} + U_{j-1} < 0$. Since
$(V^{\lambda_1} +   U_{j})(z) \to -\infty$ as $z \to
\infty$ it follows that on $\Omega_1$ we have
$(V^{\lambda_1} + U_{j})(z) < 0$. Therefore,
$(V^{\overline{\mu}_1} + \xi_{j}) < 0$ on $\Omega_1$. With
this we conclude the proof. \hfill $\Box$

\begin{rem} Notice that the functions $\widehat{s}_{1,j}(z) - \frac{P_{{\bf n},j}(z)}{Q_{\bf n}(z)}, {\bf n} \in \Lambda, $ are holomorphic in $\overline{\mathbb{C}} \setminus \Delta_1$. Using the maximum principle and \eqref{eq:6.4}, it readily follows that the sequence $({P_{{\bf n},j}}/{Q_{\bf n}}), {\bf n} \in \Lambda,$ converges to $\widehat{s}_{1,j}$   with geometric rate uniformly on any compact subset of $\mathbb{C} \setminus \Delta_1$ for each $j=1,\ldots,m$.
\end{rem}

\section{Ratio asymptotic.}

In the study of the ratio asymptotic of Hermite-Pad\'e approximants conformal mappings on Riemann surface and boundary problems of analytic functions come into play.

\subsection{Preliminaries from Riemann surfaces and boundary value problems}

Let $\Delta_1,\ldots,\Delta_m$ be a collection of intervals contained in the real line as in Definition \ref{Nikishin}. Consider the $(m+1)$-sheeted Riemann surface
$$
\mathcal R=\overline{\bigcup_{k=0}^m \mathcal R_k} ,
$$
formed by the consecutively ``glued'' sheets
$$
\mathcal R_0:=\overline {\mathbb{C}} \setminus \Delta_1,\quad
\mathcal R_k:=\overline {\mathbb{C}} \setminus \{\Delta_k \cup
\Delta_{k+1}\},\qquad k=1,\dots,m-1,\qquad \mathcal R_m=\overline
{\mathbb{C}} \setminus \Delta_m,
$$
where the upper and lower banks of the slits of two neighboring
sheets are identified. Fix $\ell \in \{1,\ldots,m\}$.  Let
$\psi^{(\ell)}, \ell=1,\ldots,m,$ be a single valued rational function
on $\mathcal{R}$ whose divisor consists of one simple zero at the
point $\infty^{(0)} \in \mathcal R_0$ and one simple pole at the
point $\infty^{(\ell)} \in \mathcal R_\ell$. Therefore,
\begin{equation}            \label{eq:psi}
\psi^{(\ell)}(z) = C_1/z + \mathcal{O}(1/z^2)\,,\,\,z \to
\infty^{(0)}\,, \qquad  \psi^{(\ell)}(z) = C_2 z +
\mathcal{O}(1)\,,\,\,z \to \infty^{(\ell)} \,,
\end{equation}
where $C_1$ and $C_2$ are constants different from zero.  Since
the genus of $\mathcal R$ equals zero (it is conformally equivalent to $\overline{\mathbb{C}}$), such a single valued
function on $\mathcal R$ exists and is uniquely determined up to a
multiplicative constant. We denote the branches of the algebraic
function $\psi^{(\ell)}$, corresponding to the different sheets $k =
0,\ldots,m$ of $\mathcal R$ by
$$
\psi^{(\ell)}:=\{\psi_k^{(\ell)}\}_{k=0}^m\,.
$$
In the sequel, we fix the multiplicative constant in such a way
that
\begin{equation} \label{eq:normaliz}
\prod_{k=0}^m |\psi_k^{(\ell)}(\infty)| = 1\,, \qquad C_1 > 0.
\end{equation}

Since $\psi^{(\ell)} $ is such that $C_1 >0,$
then
\[ \psi^{(\ell)} (z) = \overline{\psi^{(\ell)} (\overline{z})}, \qquad z \in
\mathcal{R}.
\]
In fact, define $\phi(z) := \overline{\psi^{(\ell)} (\overline{z})}$.
Notice that $\phi$ and $\psi^{(\ell)} $ have the same divisor (same poles and zeros counting multiplicities); consequently, there
exists a constant $C$ such that $\phi= C\psi^{(\ell)} $. Comparing the
leading coefficients of the Laurent expansion of these two functions
at $\infty^{(0)}$,  we conclude that $C=1$.

\medskip

In terms of the branches of $\psi^{(\ell)}$, the symmetry formula
above means that  for each  $k= 0,1,\ldots,m$:
\[
\psi_k^{(\ell)}: \overline{\mathbb{R}} \setminus
( {\Delta}_k\cup  {\Delta}_{k+1})
\longrightarrow \overline{\mathbb{R}}
\]
$( {\Delta}_0 =  {\Delta}_{m+1}=\emptyset)$;
therefore, the coefficients (in particular, the leading one) of the
Laurent expansion at $\infty$ of the branches are real numbers, and
\begin{equation}  \label{contacto}\psi_k^{(\ell)}(x_{\pm}) = \overline{\psi _k^{(\ell)}(x_{\mp})} =
\overline{\psi_{k+1}^{(\ell)}(x_{\pm})}, \qquad x \in
{\Delta}_{k+1}. \end{equation}
Among other things, the symmetry property entails
that all the coefficients in the Laurent expansion at infinity of the branches $\psi^{(\ell)}_k$ are real numbers.

\medskip

Since $\lim_{x \to \infty} x \psi_0^{(\ell)}(x) = C_1 > 0$, by continuity it follows that $\psi_k^{(\ell)}(\infty) > 0, k=1,\ldots,\ell-1, $
$\lim_{x\to \infty} \psi^{(\ell)}_\ell(x)/x = (\psi^{(\ell)}_\ell)'(\infty) > 0,$ and $\psi_k^{(\ell)}(\infty) < 0, k=\ell + 1,\ldots,m $. On the other hand, the product of all the branches $\prod_{k=0}^{m}
\psi_{k}^{(\ell)}$ is a single valued analytic function on
$\overline{\mathbb{C}}$ without singularities; therefore, by Liouville's
Theorem it is constant. Due to the previous remark and the normalization adopted in \eqref{eq:normaliz}, we can assert that
\begin{equation}\label{ident} \prod_{k=0}^{m}\,\psi_{k}^{(\ell)}(z) \equiv
\left\{
\begin{array}{rl}
1, & m-\ell\,\, \mbox{is even,} \\
-1, & m -\ell\,\, \mbox{is odd.}
\end{array}
 \right.
\end{equation}

In \cite[Lemma 4.2]{AptLopRoc} the following boundary value problem was proved to have a unique solution. In \eqref{formula} below, we introduce a slight correction  to the formula in  \cite{AptLopRoc}.

\begin{lemma}                \label{lm:system}
Let $\ell \in \{1,\ldots,m\}$ be fixed. There exists a unique collection of functions $(F_k^{(\ell)})_{k=1}^m$ which verify the system of boundary value problems
\begin{equation}          \label{eq:syst}
 \aligned &1 )\quad
F_k^{(l)},1/F_k^{(\ell)} \in \mathcal{H}(\mathbb{C} \setminus \Delta_k)\,,
\\
&2 )\quad (F _k^{(\ell)})'(\infty) > 0, \quad k=1,\dots,l\,,
\\
&2')\quad  F_k^{(\ell)}(\infty) > 0,  \quad k=l+1,\dots,m\,,
\\
&3 )\quad |F_k^{(\ell)}(x)|^2
\frac1{\bigl|(F_{k-1}^{(\ell)}F_{k+1}^{(\ell)})(x)\bigr|}= 1
 \,, \quad x \in \Delta_k\,,
\endaligned
\end{equation}
where ${F_0^{(l)}} \equiv {F_{m+1}^{(l)}} \equiv 1$. Moreover
\begin{equation} \label{formula}
F_k^{(l)} = \mbox{sg}\left(\prod_{\nu=k}^m \psi_{\nu}^{(l)}(\infty)\right)  \prod_{\nu=k}^m \psi_{\nu}^{(l)},
\end{equation}
where $\mbox{sg}\left(\prod_{\nu=k}^m \psi_{\nu}^{(l)}(\infty)\right)$ denotes the sign of the leading coefficient of the Laurent expansion at $\infty$ of
$\prod_{\nu=k}^m \psi_{\nu}^{(l)}$.
\end{lemma}

We are ready to state a result on the ratio asymptotic for type II Hermite-Pad\'e polynomials of a Nikishin system. This result was obtained in \cite{AptLopRoc} (see also \cite{LL}).

\begin{teo} \label{teo2}
Assume that $\sigma_k^{\prime} > 0$ almost everywhere on
$\Delta_k = \supp{\sigma_k}, k=1,\ldots,m$. Let $\Lambda \subset
{\mathbb{Z}}_+^m(\bullet)$ be a sequence of multi-indices such
that  $n_{1} - n_{m} \leq  d$ for all ${\bf n} \in \Lambda$, where $d$ is some fixed constant.  Then for each fixed $k \in \{1,\ldots,m\}$, we
have
\begin{equation} \label{eq:xe} \lim_{{\bf n}\in
{\Lambda}}\frac{Q_{{\bf n}^{l},k}(z)}{Q_{{\bf n},k}(z)}=
\widetilde{F_k^{(l)}}(z),
\end{equation}
uniformly on each compact subset of ${\mathbb{C}}
\setminus \Delta_k $, where $F_k^{(l)}$ is given in \eqref{formula}, the algebraic functions $\psi_{\nu}^{(l)}$ are defined by
$(\ref{eq:psi})-(\ref{eq:normaliz})$ and $\widetilde{F_k^{(l)}}$ is the result of dividing ${F_k^{(l)}}$ by the leading coefficient of its Laurent expansion at $\infty$.
\end{teo}

{\bf Sketch of the proof.} Fix $\ell \in \{1,\ldots,m\}$. Because of the interlacing property of the zeros of the polynomials $Q_{{\bf n},k}$  and $Q_{{\bf n}^\ell,k}$ it follows that for each fixed $k \in \{1,\ldots,m\}$ the family of functions $(Q_{{\bf n}^\ell,k}/Q_{{\bf n},k}), n \in \Lambda,$ is uniformly bounded on each compact subset of $\mathbb{C} \setminus \Delta_k$. To prove the theorem we need to show that for any $\Lambda' \subset \Lambda$ such that
\[ \lim_{{\bf n} \in \Lambda'} \frac{Q_{{\bf n}^\ell,k}}{Q_{{\bf n},k}} = G_k^\ell, \qquad k=1,\ldots,m,
\]
the limiting functions do not depend on $\Lambda'$.

To prove this, it is shown that their exist positive constants $c_1,\ldots,c_m$ such that $(c_k G_k^\ell)_{k=1}^m$ verifies the system of boundary value problems \eqref{eq:syst}. Properties 1), 2), and 2') are easily verified by $(G_k)_{k=1}^m$ with 1 on the right hand side of 2) and 2'). Thanks to the orthogonality properties contained in \eqref{rel3}, using results on  ratio and relative asymptotic of orthogonal polynomials with respect to varying measures contained in \cite[Theorem 6]{kn:B-G} and  \cite[Theorem 3.2]{kn:B-G2} one can also prove that $(G_k)_{k=1}^m$ satisfies $3)$ with a constant different from $1$ on the right hand side.  Normalizing the functions $G_k$ appropriately one obtains all the boundary conditions and it follows that $c_kG_k^\ell = F_k^\ell, k=1,\ldots,m$. Then using that $G_k^\ell (\infty) = 1, k=\ell +1,\ldots,m$ and $(G_k^\ell))'(\infty) = 1, k=1,\ldots,\ell$. one sees that $G_k^\ell$ has to be the function on the right hand side of \eqref{eq:xe} independently of $\Lambda'$. \hfill $\Box$

\begin{rem} An interesting open question is if one can relax the assumption $n_1 - n_m \leq d$ in the theorem. That restriction is connected with the conditions which have been found to be sufficient for the ratio and relative asymptotic of polynomials orthogonal with respect to varying measures. Those theorems in \cite{kn:B-G,kn:B-G2} would have to be improved.   Perhaps this could be done without too much difficulty if $n_1- n_m = o(|{\bf n}|), |{\bf n}| \to \infty$. Sequences verifying something like \eqref{indicesdec}, as in the case of weak asymptotic, would require a deep consideration and substantial new ideas.
\end{rem}

\begin{rem}
On the basis  of \eqref{eq:xe} and \eqref{rel4} one can also prove ratio asymptotic for the sequences $(K_{{\bf n}^\ell,k}/K_{{\bf n},k}),(\Psi_{{\bf n}^\ell,k}/\Psi_{{\bf n},k}), {\bf n} \in \Lambda, k=1,\ldots,m$ (see \cite{AptLopRoc, LL}). Using Theorem \ref{teo:4} and Proposition \ref{teo:5}, it is also possible to prove ratio asymptotic for the polynomials $A_{{\bf n},k}$ and the forms $\mathcal{A}_{{\bf n},k}$ and it is a good exercise (see also \cite{FLLS}).
\end{rem}

\begin{rem} The ratio asymptotic of multiple orthogonal polynomials finds applications in the study of asymptotic properties of modified Nikishin systems and the corresponding Hermite-Pad\'e approximants (see, for example, \cite{LL2} and \cite{LS2}).
\end{rem}

\end{document}